\documentclass[reqno,final,a4paper,11pt]{amsart}

\usepackage{amsmath,amssymb,amsthm}
\usepackage{aliascnt}
\usepackage[noadjust]{cite}
\usepackage{fullpage}
\usepackage[T1]{fontenc}
\usepackage[shortlabels]{enumitem}
\usepackage{lineno}
\usepackage{graphicx}
\usepackage[colorlinks=true,allcolors=blue]{hyperref}
\usepackage[noabbrev,capitalize,nameinlink]{cleveref}

\allowdisplaybreaks

\newcommand{\N}{\mathbb{N}}

\newtheorem{theorem}{Theorem}[section]

\newaliascnt{proposition}{theorem}
\newtheorem{proposition}[proposition]{Proposition}
\aliascntresetthe{proposition}
\crefname{proposition}{Proposition}{Propositions}

\newaliascnt{lemma}{theorem}
\newtheorem{lemma}[lemma]{Lemma}
\aliascntresetthe{lemma}
\crefname{lemma}{Lemma}{Lemmas}

\newaliascnt{claim}{theorem}
\newtheorem{claim}[claim]{Claim}
\aliascntresetthe{claim}
\crefname{claim}{Claim}{Claims}

\theoremstyle{definition}
\newaliascnt{definition}{theorem}
\newtheorem{definition}[definition]{Definition}
\aliascntresetthe{definition}
\crefname{definition}{Definition}{Definitions}

\title{Prescribed Color Patterns in Tight Hamilton Cycles and Spanning Trees}
\author{Lirong Deng}
\address{School of Mathematics and Statistics, Beijing Institute of Technology, China}
\email{lirong.deng@bit.edu.cn}
\date{}

\begin{document}
% \linenumbers
\begin{abstract}
We compare transversal and pattern thresholds for two spanning structures. 
For every fixed $k\ge 3$ and $\alpha>0$, we prove that every sufficiently large $n$-vertex $k$-graph collection with minimum codegree at least $(1/2+\alpha)n$ contains a tight Hamilton cycle with any prescribed color pattern. 
Hence the transversal and pattern thresholds for tight Hamilton cycles coincide asymptotically.
For spanning trees, we show that for every $\alpha>0$ there exists $c=c(\alpha)>0$ such that every $n$-vertex graph collection with minimum degree at least $(1/2+\alpha)n$ contains every $n$-vertex tree with maximum degree at most $n^c$ and any prescribed color pattern.  
On the other hand, we construct spanning trees of maximum degree $\sqrt n+1$ whose pattern threshold exceeds their transversal threshold by a positive linear term.
\end{abstract}

\maketitle
\section{Introduction}

Minimum degree conditions forcing spanning structures are a central topic in extremal graph theory.
A classical example is Dirac’s theorem~\cite{MR47308}, which states that every graph on $n\ge3$ vertices with minimum degree at least $n/2$ contains a Hamilton cycle. For uniform
hypergraphs, the corresponding questions depend both on the notion of a cycle and on the degree parameter under consideration.

For a $k$-uniform hypergraph ($k$-graph for short) $H$ and $U\subseteq V(H)$, let $H[U]$ denote the $k$-graph induced by $U$.
For $S\in \binom{V(H)}{k-1}$, define $N_H(S):=\{v\in V(H)\setminus S:S\cup \{v\}\in E(H)\}$ and $d_H(S):=|N_H(S)|$.
The \emph{minimum codegree} of $H$ is $\delta_{k-1}(H)= \min_{S\in\binom{V(H)}{k-1}}d_H(S)$.
A \emph{tight Hamilton cycle} in $H$ is a cyclic ordering of $V(H)$ such that every set of $k$ consecutive vertices forms an edge of $H$.
R{\"o}dl, Ruci{\'n}ski and Szemer{\'e}di~\cite{rodl2008approximate} proved that, for every fixed $k\ge3$ and $\alpha>0$, every sufficiently large $n$-vertex $k$-graph with minimum codegree at least $(1/2+\alpha)n$ contains a tight Hamilton cycle. 
The coefficient $1/2$ is asymptotically best possible.

A natural extension is to replace a single host graph by a collection of graphs on the same vertex set.
Let $\mathcal H=\{H_1,\ldots,H_m\}$ be a collection of $k$-graphs on a common vertex set $V$, where repeated members are allowed, and write $\delta_{k-1}(\mathcal{H}):=\min_{i\in [m]}\delta_{k-1}(H_i)$. 
For a positive integer $m$, write $[m]:=\{1,\dots,m\}$. 
For $R\subseteq V$, we further define $\delta_{k-1}(\mathcal{H},R):=\min_{i\in [m]\ S\in\binom{V}{k-1}}|N_{H_i}(S)\cap R|$.
For $t\in[m]$, let $\mathcal{H}[t]:=(H_1,\dots,H_t)$ when we consider the order of these $k$-graphs.

\subsection{Transversal thresholds}
A copy $F'$ of a $k$-graph $F$ is \emph{transversal}, or \emph{rainbow}, in $\mathcal{H}$ if there is an injection $\varphi: E(F')\to [m]$ such that $e\in E(H_{\varphi(e)})$ for every $e\in E(F')$.
When $m=e(F)$, every member of the system is used exactly once. Taking all members to be identical recovers the ordinary embedding problem.
Consequently, a necessary degree condition for a transversal embedding is the corresponding condition for a single host. Whether this necessary condition is also sufficient is the question of \emph{color-blindness}.

A fundamental result of Joos and Kim~\cite{joos2020rainbow} establishes exact color-blindness for Hamilton cycles in graphs: every collection of $n$ graphs on a common $n$-vertex set with minimum degree at least $n/2$ contains a transversal Hamilton cycle.
Montgomery, M\"uyesser and Pehova~\cite{MR4451150} proved that, for every $\alpha>0$, there exists $c>0$ such that every collection of $n-1$ graphs with minimum degree at least $(1/2+\alpha)n$ contains a transversal copy of every $n$-vertex tree of maximum degree at most $cn/\log n$, provided that $n$ is sufficiently large. 
Their work also treats factors of fixed graphs and shows that color-blindness does not hold for every such factor.

For hypergraphs, Cheng, Han, Wang, Wang and Yang~\cite{cheng2025transversal} showed that minimum codegree $(1/2+\alpha)n$ in each member of a collection of $n$ $k$-graphs guarantees a transversal tight Hamilton cycle, for fixed $k\ge3$, $\alpha>0$ and sufficiently large $n$.
Thus tight Hamilton cycles are asymptotically color-blind with respect to minimum codegree.
Transversal versions of other spanning embedding results have also been obtained, including clique factors and hypergraph perfect matchings~\cite{cheng2023rainbow}, powers of Hamilton cycles~\cite{gupta2023general}, graphs of small bandwidth~\cite{chakraborti2023bandwidth}, and Hamilton
paths and cycles in tournaments~\cite{chakraborti2024hamilton}.

\subsection{Pattern thresholds}

The pattern problem strengthens the transversal problem by fixing the assignment of the members of the collection to the edges of the desired spanning structure in advance.
In the language of Bowtell, Morris, Pehova and Staden~\cite{bowtell2025universality}, one asks not merely for some transversal copy, but for each prescribed color pattern, a copy realising that pattern.
For a tight Hamilton cycle $C=e_1\cdots e_n$ in a collection $\mathcal H=\{H_1,\ldots,H_m\}$, the pattern is naturally encoded by a map $\chi:[n]\to[m]$, and the requirement is $e_i\in E(H_{\chi(i)})$ for every $i\in[n]$.
When $\chi$ is injective, this is stronger than the transversal problem, because the member assigned to each position of the cycle is prescribed rather than chosen by the embedding.

For Hamilton cycles in graphs, Bowtell et al.~\cite{bowtell2025universality} proved that, for every fixed $\alpha>0$, minimum degree at least $(1/2+\alpha)n$ in every member of the collection guarantees every prescribed Hamilton cycle pattern. 
Hence the classical, transversal and pattern minimum-degree thresholds for Hamilton cycles all have leading term $n/2$.
Their result also shows why the pattern problem is genuinely stronger: at the exact Dirac threshold $n/2$, the construction of Gupta et al., reported in [\cite{bowtell2025universality}, Section 1.4], excludes a prescribed Hamilton-cycle pattern.
Taking two copies of its first graph and $n-2$ copies of its second gives a
collection admitting a transversal Hamilton cycle by Joos-Kim, but excluding the ordering in which the two first-graph occurrences are consecutive.
Thus coincidence of the leading terms need not imply equality of the exact thresholds.

The comparison between transversal and pattern thresholds has already appeared in the study of factors. 
In particular, Montgomery et al.~\cite{MR4451150} gave examples where the two thresholds coincide and examples where they are separated. 
More recently, universality questions for prescribed patterns have been considered for powers of Hamilton cycles and in random graph collections; see, for example,~\cite{heath2025universality,christoph2025universality}. These results indicate that the relation between transversal and pattern thresholds can depend strongly on the local structure of the target and on how much freedom an embedding has to rearrange its colors.

\subsection{Main results}

In the concluding remarks of ~\cite[Section~7.2]{bowtell2025universality}, Bowtell et al.~highlighted Hamilton cycles in hypergraphs and bounded-degree spanning trees as natural directions for extending pattern universality.
Motivated by these directions, we investigate the relation between transversal and pattern thresholds for both spanning structures.

Our first result extends their graph theorem to tight Hamilton cycles in uniform hypergraphs, showing that prescribing the complete color pattern
causes no asymptotic cost under the Dirac codegree condition.

\begin{theorem}\label{thm:pattern Hamilton}
For every $\alpha>0$ and every integer $k\ge 3$, there exists $n_0$ such that for every $n\ge n_0$ the following holds. 
Let $m\in\N$ and let $\mathcal{H}=\{H_1,\dots,H_m\}$ be a collection of $k$-graphs on a common vertex set $V=[n]$, and suppose that $\delta_{k-1}(\mathcal{H})\ge (1/2+\alpha)n$.
Then, for every $\chi$: $[n]\to [m]$, there exists a tight Hamilton cycle $C=e_1\cdots e_n$ such that $e_i\in H_{\chi(i)}$ for every $i\in [n]$.
\end{theorem}

We use round brackets to emphasize the implicit ordering of $\{H_1,\dots,H_m\}$ which encodes this canonical coloring.
As in~\cite{bowtell2025universality}, repeated members allow us to reduce ~\cref{thm:pattern Hamilton} to the canonical ordered pattern.
More precisely, it is enough to prove the following statement.

\begin{theorem}\label{thm:sufficient main result}
    For every $\alpha>0$ and every integer $k\ge 3$, there exists $n_0$ such that for every $n\ge n_0$ the following holds. 
    Let $\mathcal{H}=(H_1,\dots,H_n)$ be an ordered  collection of $k$-graphs on a common vertex set $V=[n]$, and suppose that $\delta_{k-1}(\mathcal{H})\ge (1/2+\alpha)n$.
    Then, there exists a tight Hamilton cycle $C=e_1\cdots e_n$ such that $e_i\in H_i$ for all $i\in [n]$.
\end{theorem}

\Cref{thm:pattern Hamilton} and~\Cref{thm:sufficient main result} show that, under the asymptotic Dirac codegree condition, the pattern threshold of tight Hamilton cycles meets the transversal threshold.

The same asymptotic coincidence holds for spanning trees of sufficiently small polynomial maximum degree.

\begin{theorem}\label{thm:tree}
For every $\alpha>0$, there exist $c=c(\alpha)>0$ and $n_0$ such that the following holds for all $n\ge n_0$.
Let $\mathcal G=(G_1,\ldots,G_m)$ be a collection of graphs on an $n$-vertex set $V$ satisfying $\delta(\mathcal G)\ge (1/2+\alpha)n.$
Then, for every $n$-vertex tree $T$ with $\Delta(T)\le n^c$ and every map $\chi:E(T)\to[m]$, there exists a bijection
$f:V(T)\to V$ such that $f(x)f(y)\in E(G_{\chi(xy)})$ for every $xy\in E(T)$.
\end{theorem}

Our next result shows that the coincidence in ~\cref{thm:tree} does not extend to all trees of sublinear maximum degree.
More precisely, we construct a sequence of spanning trees of sublinear maximum degree for which a prescribed pattern requires a linearly larger minimum degree than a transversal copy.

\begin{theorem}\label{thm:tree-construction}
For every integer $\ell\ge 9$, there exist arbitrarily large square
integers $n$ such that the following holds.
Let $m=3\ell+1$.
There exists a collection $\mathcal G=\{G_1,\dots, G_m\}$ of graphs on a common vertex set $V=[n]$, an $n$-vertex tree
$T$, and a map $\chi:E(T)\to [m]$ such that
$\delta(\mathcal G)\ge\left(\frac12+\frac{1}{2(24\ell+1)}\right)n-1$
and $\Delta(T)=\sqrt n+1$,
but there is no bijection $f:V(T)\to V$ such that $f(x)f(y)\in E(G_{\chi(xy)})$ for every $xy\in E(T)$.
\end{theorem}

The significance of~\Cref{thm:tree-construction} is highlighted by comparison with the transversal theorem of Montgomery et al.~\cite{MR4451150}. 
Replacing each prescribed occurrence of a color on an edge $e\in E(T)$ by a separate copy of $G_{\chi(e)}$ gives a collection of $n-1$ graphs with the same minimum degree. 
Since $\Delta(T)=\sqrt n+1=o(n/\log n)$, this collection contains a transversal copy of $T$ for every fixed positive margin above $n/2$ and all sufficiently large $n$, whereas~\Cref{thm:tree-construction} rules out the prescribed pattern. 
In particular, for $\ell=9$, the pattern threshold is at least $(1/2+1/434)n-O(1)$ along the corresponding sequence of orders, while the transversal threshold is $(1/2+o(1))n$. 
Thus, within the two settings considered here, we see two different behaviors. 
For tight Hamilton cycles, sufficient local flexibility allows the prescribed color order to be maintained at no asymptotic cost, whereas for the trees in ~\Cref{thm:tree-construction}, local restrictions imposed by the prescribed colors accumulate into a global obstruction.

\subsection{Organisation}

In ~\Cref{sec:notation}, we introduce the notation and preliminary results. 
In ~\Cref{sec:proof}, we state the main lemmas and deduce ~\Cref{thm:sufficient main result} and ~\cref{thm:tree}. 
~\Cref{sec:connecting,sec:absorption,sec:almost path cover} prove the connecting, absorbing and almost-cover lemmas, respectively.
~\cref{sec:tree} contains the two additional lemmas needed for spanning trees.
~\Cref{sec:construction} contains the construction proving ~\Cref{thm:tree-construction}. 
We conclude in ~\Cref{sec:conluding rmk} with several questions suggested by the comparison between transversal and pattern thresholds.

\section{Notation and preliminaries}\label{sec:notation}

Throughout the paper, floors and ceilings are omitted in the interest of readability. 
In any such case, the result can be made precise without non-trivial changes to the proof.
We write $x\ll y\ll z$ to mean that we choose constants hierarchically from right to left, that is, for any $z>0$, there exist functions $f$ and $g$ such that, whenever $y\le f(z)$ and $x\le g(y)$, the subsequent statement holds. 
Hierarchies of other lengths are defined analogously.

A $(k-1)$-\emph{tuple} is a $(k-1)$-element sequence of distinct vertices of $V$.
For a finite set $U$ and an integer $t\ge0$, let $\binom{U}{t}$ denote the family of $t$-subsets of $U$, and let $\Omega_t(U):=\{(u_1,\dots,u_t)\in U^t:u_1,\dots, u_t\ \text{are distinct}\}$
denote the set of ordered $t$-tuples of distinct elements of $U$.
In particular, we set $\Omega_0(U):=\{()\}$, where $()$ denotes the empty tuple.

Let $H$ be a $k$-graph.
A \emph{tight path} $P=v_1v_2\cdots v_t$ in $H$ is an ordered sequence of distinct vertices such that $e_i(P)=\{v_i,v_{i+1},\dots, v_{i+k-1}\}\in E(H)$ for every $i\in [t-k+1]$.
The \emph{length} of $P$ is its number of edges.
The ordered $(k-1)$-tuples $\operatorname{in}(P):=(v_1,\dots,v_{k-1})$ and $\operatorname{ter}(P):=(v_{t-k+2},\dots,v_t)$ are called the \emph{initial} and \emph{terminal ends} of $P$, respectively.  
For ordered $(k-1)$-tuples $\mathbf{x}$ and $\mathbf{y}$, an $\mathbf{x},\mathbf{y}$-tight path is a tight path $P$ satisfying $\operatorname{in}(P)=\mathbf{x}$ and $\operatorname{ter}(P)=\mathbf{y}$.

Let $\mathcal{H}=(H_1,\dots,H_m)$ be an ordered collection of $k$-graphs on a common vertex set.
Suppose that $P$ has ordered edges $e_1,\dots,e_{\ell}$ and let $\psi: [\ell]\to [m]$.
We say that $P$ is \emph{exactly $\psi$-colored} in $\mathcal{H}$ if $e_j\in E(H_{\psi(j)})$ for every $j\in [\ell]$.
When $m=\ell$ and $\psi$ is the identity map, we simply say that $P$ is \emph{exactly $\mathcal{H}$-colored}.  

For an ordered $(k-1)$-tuple $\mathbf{x}=(x_1,\dots,x_{k-1})\in\Omega_{k-1}(V)$ and a vertex $y\notin \{x_1,\dots,x_{k-1}\}$, write $\mathbf{x}\circ y:=(x_2,\ldots,x_{k-1},y).$
If $\{x_1,\ldots,x_{k-1},y\}\in E(H)$, then $y$ is called a \emph{legal extension} of $\mathbf x$ in $H$.  

Let $H$ be a $k$-graph.
A function $\omega: E(H)\to \mathbb{R}^+$ is a \emph{fractional matching} of $H$ if $\sum_{e \ni v} \omega(e) \le 1$ for every $v\in V(H)$. 
It is \emph{perfect} if equality holds for every $v\in V(H)$.
If $|V(H)|=M$, we call a perfect fractional matching $\omega$ \emph{$C$-normal} if $\frac{1}{CM^{k-1}}\le \omega(e)\le \frac{C}{M^{k-1}}$ for every $e\in E(H)$.

The following result will allow us to choose extensions of an ordered $(k-1)$-tuple in a balanced way. 
This will be the basic ingredient in the random walks used later.
For this purpose, we shall use the following consequence of Cuckler and Kahn ~\cite{cuckler2009hamiltonian} on entropy-maximising perfect fractional matchings.
The stated formulation is the graph case of normal perfect fractional matching lemma of Glock, Gould, Joos, K\"uhn and Osthus~\cite{glock2021counting}.

\begin{lemma}[Graph case of {\cite[Lemma 4.2]{glock2021counting}}]\label{lem:normal-fpm}
Let $1/n\ll1/C\ll\alpha$, and let $G$ be a graph on $n$ vertices with $\delta(G)\ge \left(\frac12+\alpha\right)n$.
Then there exists a $C$-normal perfect fractional matching of $G$.
\end{lemma}

We will also need two standard concentration inequalities. 
The first is the usual Chernoff bound for binomial and hypergeometric random variables.

\begin{lemma}[Chernoff bounds {\cite[Corollary 21.7]{frieze2015introduction}}]\label{lem:chernoff}
    Let $0<\varepsilon<1$ and let $X$ be a binomially or hypergeometrically distributed random variable with mean $\mu$. Then
    \[\Pr(|X-\mu|\ge \varepsilon \mu)\le 2e^{-\varepsilon^2\mu/3}.\]
\end{lemma}

A sequence $(X_i)_{i=0}^n$ of random variables is a \emph{martingale} if $\mathbb{E}[X_{i+1}\mid X_i,X_{i-1},\dots,X_0]=X_{i}$ for every $0\le i<n$.
We shall use the standard Azuma--Hoeffding inequality; see, for example, \cite[Theorem 7.2.1 and its proof]{AlonSpencer2016}.

\begin{lemma}[Azuma--Hoeffding inequality]\label{lem:Azuma}
Let $(X_i)_{i=0}^n$ be a martingale with respect to a filtration
$(\mathcal F_i)_{i=0}^n$. 
Suppose that, for every $i\in[n]$, $|X_i-X_{i-1}|\le c_i$
almost surely, where $c_1,\ldots,c_n$ are positive constants.
Then, for every $t>0$,
\[
    \Pr\bigl(|X_n-X_0|\ge t\bigr)
    \le
    2\exp\left(
        -\frac{t^2}{2\sum_{i=1}^n c_i^2}
    \right).
\]
\end{lemma}

Finally we need the following result which shows that we can select a subset of the vertex set of a hypergraph that maintains the minimum codegree condition.
In ~\cite{bowtell2025universality}, this result was established for graphs. 
Here we state its straightforward extension to hypergraphs.

\begin{proposition}[{\cite[Proposition 2.3]{bowtell2025universality}}]\label{prop:subset}
Let $k\ge 2$ be a fixed integer and let
$0<1/n\ll \alpha,\beta,1/k<1$.
Let $\mathcal{H}=(H_1,\ldots,H_n)$ be a collection
of $k$-graphs on a common vertex set $V=[n]$ with $\delta_{k-1}(\mathcal{H})\ge (1/2+\alpha)n.$
Then there exists $Z\subseteq V$ such that
$|Z|=\lfloor\beta n\rfloor$ and
\[
\delta_{k-1}(\mathcal{H},Z)
\ge (1/2+\alpha/2)|Z|.
\]
\end{proposition}
\begin{proof}
    Let $z=\lfloor\beta n\rfloor$ and choose a set
    $Z\in\binom{V}{z}$ uniformly at random.
    For each $i\in[n]$ and $S\in\binom{V}{k-1}$, the random variable $X_{i,S}:=|N_{H_i}(S)\cap Z|$ has a hypergeometric distribution with mean $\mu_{i,S}=zd_{H_i}(S)/n \ge \left(\frac12+\alpha\right)z$.
    By \cref{lem:chernoff}, taking $\varepsilon=\alpha z/(2\mu_{i,S})$, we obtain $\Pr\left(X_{i,S}<\left(\frac12+\frac{\alpha}{2}\right)z\right)\le
    2\exp\left(-\frac{\alpha^2z^2}{12\mu_{i,S}}\right)
    \le 2\exp\left(-\frac{\alpha^2z}{12}\right)$,
    where the last inequality follows from
    $\mu_{i,S}\le z$.
    Since there are at most
    $n\binom{n}{k-1}\le n^k$ choices of $(i,S)$,
    a union bound gives
    \[
    \Pr\left(
    \delta_{k-1}(\mathcal H,Z)
    <\left(\frac12+\frac{\alpha}{2}\right)z
    \right)
    \le
    2n^k\exp\left(-\frac{\alpha^2z}{12}\right)
    =o(1).
    \]
    Thus the desired set $Z$ exists.
\end{proof}

\section{Proof of the main results}\label{sec:proof}

In this section, we state the main ingredients of the proof, explain how they fit together, and then deduce ~\cref{thm:sufficient main result} and ~\cref{thm:tree} from them. 
Their proofs are deferred to ~\cref{sec:connecting,sec:absorption,sec:almost path cover,sec:tree}.

\subsection{Main lemmas for tight Hamilton cycles}

We first include a connecting lemma, which allows us to join two prescribed ordered ends while respecting an arbitrary prescribed interval of colors.

\begin{lemma}[Connecting]\label{lem:connecting}
    Let $0<1/N\ll1/r\ll\alpha, 1/k$, where $k\ge2$. 
    Let $R$ be an $N$-set, $R'\subseteq R$ with $|R\setminus R'|\le \frac{\alpha}{2}N$. 
    Let $\mathcal{J}=(J_1,\dots,J_{r+k-1})$ be an ordered collection of $k$-graphs on a vertex set containing $R$ with $\delta_{k-1}(\mathcal{J},R)\ge (1/2+\alpha)N$.
    Then any two disjoint ordered $(k-1)$-tuples are joined by $\Omega(|R'|^r)$ exactly $\mathcal{J}$-colored tight paths whose internal vertices lie in $R'$. 
\end{lemma}

Our second ingredient is an absorbing path whose color pattern is fixed in advance, while its vertex set retains sufficient flexibility to absorb a prescribed subset of a small reservoir.

\begin{lemma}[Absorber]\label{lem:absorber}
    Let $0<1/n\ll\beta\ll\alpha, 1/k$, where $k\ge2$.
    Then there exists an integer $a\le\alpha n/20$ such that the following holds.
    Let $\mathcal{J}=(J_1,\dots,J_a)$ be an ordered collection of $k$-graphs on a vertex set $V=[n]$ with $\delta_{k-1}(\mathcal{J})\ge (1/2+\alpha)n$.
    For every $Z\subseteq V$ with $\lfloor\beta n\rfloor\le |Z|\le 2\lfloor\beta n\rfloor$, there exist $A\subseteq V\setminus Z$ and disjoint ordered $(k-1)$-tuples $\mathbf{x},\mathbf{y}\in \Omega_{k-1}(A)$ with $|A|=a+k-1-\lfloor\beta n\rfloor$ such that, for every $Z'\in \binom{Z}{\lfloor\beta n\rfloor}$, there is an exactly $\mathcal{J}$-colored $\mathbf{x},\mathbf{y}$-tight path with vertex set $A\cup Z'$.
\end{lemma}

Finally, we need an almost-cover result which produces many vertex-disjoint tight paths with independently prescribed color patterns.

\begin{lemma}[Almost path cover]\label{lem:almost cover}
    Let $0<1/n\ll \varepsilon, 1/K\ll \alpha, 1/k$, where $k\ge 2$ and let $m\le n$.
    Let $\mathcal{J}=(J_1,\dots,J_m)$ be a collection of $k$-graphs on a vertex set $V=[n]$ with $\delta_{k-1}(\mathcal{J})\ge (1/2+\alpha)n$.
    Set $s:=\lfloor (1-\varepsilon)n/K\rfloor$.
    Then for any maps $\psi_j$: $[K-k+1]\to [m]$, $j\in[s]$, there exist pairwise vertex-disjoint $K$-vertex tight paths $P_1,\dots,P_s$ such that $P_j$ is exactly $\psi_j$-colored in $\mathcal{J}$ for every $j\in[s]$.
\end{lemma}

\subsection{Main lemmas for spanning trees}
We need two further lemmas.
The first deals with trees having linearly many leaves and a small polynomial maximum-degree bound. 
The second fills a collection of bare paths after their endpoints have already been embedded.

\begin{lemma}[Trees with many leaves]\label{lem:leaves}
For every $\alpha>0$ and $0<\lambda<1$, there exist $c=c(\alpha,\lambda)>0$ and $n_0$ such that the following holds
for all $n\ge n_0$.
Let $\mathcal G=(G_1,\ldots,G_m)$ be a graph collection on an $n$-vertex set $V$ with $\delta(\mathcal G)\ge(1/2+\alpha)n.$
If $T$ is an $n$-vertex tree satisfying $\Delta(T)\le n^c$ and having at least $\lambda n$ leaves, then $\mathcal G$ contains an exactly $\chi$-colored copy of $T$ for every $\chi:E(T)\to[m]$.
\end{lemma}

\begin{lemma}[Filling bare paths]\label{lem:linkage}
Let $0<1/N\ll1/K\ll\alpha$. 
Let $U$ be an $N$-set and let $X$ be disjoint from $U$.
Suppose that $\mathcal G=(G_1,\ldots,G_m)$ is a graph collection on $U\cup X$ with $\delta(\mathcal{G},U)\ge(1/2+\alpha)N$.
For $j\in[p]$, let $x_j,y_j\in X$ be distinct and let
$\ell_j\in[K,2K]$ be integers satisfying $\sum_{j=1}^p\ell_j=N.$
Then for any maps $\psi_j:[\ell_j+1]\to[m]$, there are exactly $\psi_j$-colored $x_j,y_j$-paths $P_j$ whose internal vertex sets partition $U$, with $P_j$ having precisely $\ell_j$ internal vertices.
\end{lemma}

Here different paths may share endpoints.

\subsection{Proof outline}

We briefly describe how these main lemmas are combined to prove ~\cref{thm:sufficient main result} and ~\cref{thm:tree}. 
The proof follows the usual absorption scheme, with the additional requirement that all edges respect the prescribed order of the colors. The main point is therefore to reserve suitable intervals of colors for the absorber, the almost-cover paths, and the connecting paths simultaneously.

For tight Hamilton cycles, we first choose a small reservoir $Z$ which inherits the
minimum codegree condition in every color.
Applying ~\cref{lem:absorber} to an initial interval of colors, we obtain an absorbing set $A$ with ordered ends $\mathbf{x},\mathbf{y}$. 
Its key property is that, for any prescribed $\beta n$-subset $Z'\subseteq Z$, the set $A\cup Z'$ spans an exactly colored tight $\mathbf{x},\mathbf{y}$-path.

We then work in $W:=V\setminus(A\cup Z)$. 
~\cref{lem:almost cover} gives a collection of pairwise vertex-disjoint tight paths covering almost all of $W$, with
the required color patterns. 
Using fresh vertices of $Z$, ~\cref{lem:connecting} allows us to connect these paths in the prescribed order. 
The few remaining vertices of $W$ are incorporated in the same way, after which we extend inside $Z$ and finally connect back to $\mathbf{x}$, leaving exactly $\beta n$ vertices of $Z$ unused.

For spanning trees, we distinguish two cases according to the number of leaves.

If the tree has linearly many leaves, we first remove a linear number of them and randomly embed the remaining tree in connected pieces.
Balanced extensions(Definition see~\cref{sec:connecting}) make most vertex images approximately uniform. 
This preserves the color neighborhoods and distributes the parents of the removed leaves, allowing us to embed the leaves by a perfect matching.

If the tree has few leaves, almost all its vertices lie on long bare paths. We first embed a small forest containing the endpoints of these paths, and then apply ~\cref{lem:linkage} to fill their interiors with the prescribed colors.

\subsection{Proof of \cref{thm:sufficient main result}}
\begin{proof}[Proof of \cref{thm:sufficient main result}]
    Choose $0<1/n\ll \varepsilon, 1/K\ll \gamma\ll\beta\ll\alpha,1/k$.
    By~\cref{prop:subset}, we obtain $Z\subseteq V$ of size $(\gamma+\beta)n$ such that $\delta_{k-1}(\mathcal{H},Z)\ge(1/2+\alpha/2)|Z|$.

    Let $a$ be given by \cref{lem:absorber}. 
    Applying ~\cref{lem:absorber} to the ordered
    collection $\mathcal{H}[a]=(H_1,\dots,H_a)$ and the set $Z$, we obtain $A\subseteq V\setminus Z$ and disjoint ordered $(k-1)$-tuples $\mathbf{x},\mathbf{y}\in \Omega_{k-1}(A)$, where $|A|=a+k-1-\beta n$, $a\le\alpha n/20$,
    and for every $Z'\in\binom{Z}{\beta n}$, there is an exactly
    $\mathcal{H}[a]$-colored tight $\mathbf x,\mathbf y$-path on $A\cup Z'$.
    
    Set $W:=V\setminus(A\cup Z)$ and $n':=|W|=n-a-\gamma n-k+1$.
    By the choice of constants, deleting $A\cup Z$ gives $\delta_{k-1}(H_i[W])\ge (1/2+\alpha/2)n'$ for every $i\in[n]$.

    Set $s=\lfloor(1-\varepsilon)n'/K\rfloor$, choose a number of internal vertices $r$ from ~\cref{lem:connecting} such that $r/K+\varepsilon(r+k-2)<\gamma/4$.
    For $j\in[s]$ and $t\in [K-k+1]$, define the $K-k+1$ colors of the $j$th path by $\psi_{j}(t):=a+(j-1)(K+r)+(r+k-1)+t$.
    Thus, in the $j$th block of $K+r$ consecutive colors, the first
    $r+k-1$ colors are reserved for a connector and the remaining $K-k+1$ colors are assigned to the $K$-vertex path.
    Let $I:=\{\psi_j(t):j\in[s],\ t\in[K-k+1]\}.$
    Since $|I|=s(K-k+1)\le sK\le n',$
    we may apply ~\cref{lem:almost cover} to the ordered subcollection
    $(H_i[W])_{i\in I}$, with the indices relabelled in the natural way. 
    We obtain pairwise disjoint $K$-vertex tight paths $P_1,\dots, P_s\subseteq W$ such that $P_j$ is exactly $\psi_j$-colored for every $j\in[s]$.
    Using unused vertices of $Z$, apply ~\cref{lem:connecting} to connect $\mathbf{y}\rightarrow P_1\rightarrow P_2\rightarrow\cdots\rightarrow P_s$, where each connector has $r$ internal vertices and uses the reserved interval of $r+k-1$ colors.
    Indeed, throughout the construction only $O(\gamma n)$ vertices of $Z$ are forbidden, whereas $|Z|=(\beta+\gamma)n$ and $\gamma\ll\beta,\alpha$; hence the available subset of $Z$ satisfies
    the condition in \cref{lem:connecting}.
    Denote the resulting tight path by $P$. It starts with $\mathbf y$, covers all vertices of $P_1,\ldots,P_s$, and uses precisely the colors $a+1,\ldots,a+s(K+r).$
    At this point exactly $sr$ vertices of $Z$ have been used.
    
    We next deal with the remaining $d:=n'-sK$ vertices of $W$ that have not yet been covered. 
    Write these vertices as $v_1,\dots,v_d$.
    For each $v_j$, choose $k-2$ unused vertices $z_{j,1},\dots,z_{j,k-2}\in Z$ and form the ordered tuple $(v_j,z_{j,1},\dots,z_{j,k-2})$.
    Apply ~\cref{lem:connecting} to connect the terminal end of the current tight path to this tuple using the prescribed connector with $r$ internal vertices.
    For each step, we use one vertex of $W$ and $r+k-2$ vertices of $Z$.
    After all $d$ uncovered vertices have been inserted, the number of used vertices in $Z$ is $z_0:=sr+d(r+k-2)<\gamma n/2$.

    For the remaining vertices in $Z$, greedily extend the current tight path by $p:=\gamma n-z_0-r$ unused vertices of $Z$ and always use the next prescribed color.
    $\delta_{k-1}(\mathcal{H},Z)$ guarantees a legal extension at every step.
    This leaves $\beta n+r$ unused vertices of $Z$.
    Apply~\Cref{lem:connecting} to all these unused vertices, to choose $r$ internal vertices of the final connector to reach $\mathbf{x}$.

    Set $Z':=Z\setminus V(P)$, then $|Z'|=\beta n$.
    After color $a$, the number of edges used is $s(K+r)+d(r+k-1)+(\gamma n-z_0-r)+(r+k-1)=sK+d+\gamma n+k-1=n-a$.
    Therefore the tight path from $\mathbf{y}$ back to $\mathbf{x}$ uses exactly the colors $a+1,\dots,n.$
    Finally the absorbing property of $A$ gives an exactly $\mathcal{H}[a]$-colored tight path from $\mathbf{x}$ to $\mathbf{y}$ with vertex set $A\cup Z'$.
    Concatenating the two paths cyclically yields a tight Hamilton cycle whose $i$th edge lies in $H_i$ for every $i\in[n]$.
    This proves the theorem.
\end{proof}

\subsection{Proof of ~\cref{thm:tree}}
We use the standard observation that a tree with few leaves can be covered almost entirely by long bare paths; see, for example,~\cite[Lemma~4.1]{montgomery2020embedding}.
For completeness, we give the decomposition needed here.
\begin{proof}[Proof of ~\cref{thm:tree}]
It suffices to consider $0<\alpha<1/4$. 
Choose $1/n\ll 1/K\ll\alpha$ and then choose $\lambda>0$ sufficiently small that $2(K+1)\lambda+1/(K+1)<\alpha/4$.

Let $c=c(\alpha,\lambda)$ be given by ~\cref{lem:leaves}.  
Fix a tree $T$ with $\Delta(T)\le n^c$ and a pattern $\chi:E(T)\to[m]$. 
If $T$ has at least $\lambda n$ leaves, apply ~\cref{lem:leaves} and we obtain the desired result. 
We may therefore assume that its number of leaves $L_T <\lambda n$.

Put $B_T=\{v\in V(T):d_T(v)\ne2\}$. 
By the Handshaking Lemma we have $\sum_{v\in V(T)}d_T(v)=2(n-1)$.
Hence, the identity $\sum_{v:d_T(v)\ge3}(d_T(v)-2)=L_T-2$ shows that $|B_T|\le2L_T$. 
Note that every vertex outside $B_T$ lies on a chain between two vertices of $B_T$. 
Suppressing all degree-two vertices gives a tree on $B_T$, so there are $|B_T|-1$ such chains.

We cut each chain into consecutive segments of $K+1$ edges and a final segment of at most $K$ edges. 
Let $C_T$ be the set consisting of $B_T$, the cut vertices, and all vertices in the final short segments. 
The vertices outside $C_T$ are partitioned into interiors of bare paths, each with exactly $K$ internal vertices and two distinct endpoints in $C_T$. 
There are at most $n/(K+1)$ long segments, and each chain contributes at most $K$ additional vertices from its short segment. 
Consequently,
\[
|C_T|\le |B_T|+\frac n{K+1}+K|B_T|
\le2(K+1)L_T+\frac n{K+1}<\frac{\alpha n}{4}.
\]

We embed the forest $T[C_T]$ greedily, processing each component from a root and respecting $\chi$. 
Every new vertex with an already embedded parent has at least $(1/2+\alpha)n$ possible neighbors in its prescribed color, whereas fewer than $|C_T|$ host vertices have been used. 
The root of a new component can be placed at any unused vertex. Thus this produces an embedding $f$.

Put $U=V\setminus f(C_T)$, $N=|U|$, and $X=f(C_T)$. 
For every $v\in U\cup X$ and every color $i$, $d_{G_i}(v,U)\ge(1/2+\alpha)n-|C_T|\ge(1/2+\alpha/2)N.$
All unembedded bare paths have their endpoint images fixed in $X$, their internal vertex counts are $K$, and their sum is $N$. 
Apply ~\cref{lem:linkage} to fill them with their original color lists using every vertex of $U$ exactly once. 
Together with $f$, this gives the required spanning copy of $T$.
\end{proof}

\section{Connecting}\label{sec:connecting}

In this section, we prove ~\cref{lem:connecting}. 
The main difficulty in constructing paths with prescribed color sequences is to control the distribution of their ordered ends.
Choosing each new vertex uniformly from the available neighborhood does not preserve the uniform distribution in general, particularly when the underlying hypergraph changes from one step to the next. 
To overcome this difficulty, we use normal perfect fractional matchings in the link graphs to define a balanced extension rule. 
This rule preserves the uniform distribution on ordered $(k-1)$-tuples, while assigning probabilities of $O(1/n)$ to all legal extensions.
Moreover, the minimum codegree condition ensures that two such walks can be coupled to agree after a bounded number of steps with positive probability. 
These properties allow us to connect arbitrary ordered ends through a prescribed sequence of colors, with polynomially many choices for the internal vertices.
Our construction is closely related to the weighted random-walk
approach of Glock, Gould, Joos, K\"uhn and Osthus
\cite{glock2021counting}, which uses normal perfect fractional
matchings to control the distribution of a tight walk.
Here we apply their fractional matching lemma separately in
each $(k-2)$-link graph.
We begin by defining the extension rule.

\subsection{Balanced extensions}
Let $0<1/n\ll \alpha,1/k$ and let $H$ be an $n$-vertex $k$-graph satisfying $\delta_{k-1}(H)\ge \left(1/2+\alpha\right)n$. 
Fix an ordered $(k-2)$-tuple $\mathbf s=(s_1,\ldots,s_{k-2})\in\Omega_{k-2}(V(H))$. 
Define the \emph{link graph} $L_H(\mathbf s)$ on the vertex set
$V(H)\setminus\{s_1,\ldots,s_{k-2}\}$ by $xy\in E(L_H(\mathbf s))$ if and only if $\{x,s_1,\ldots,s_{k-2},y\}\in E(H).$
Put $m:=n-k+2$.  
For every $x\in V(L_H(\mathbf s))$,
we have $d_{L_H(\mathbf s)}(x)=d_H(\{x,s_1,\ldots,s_{k-2}\})
\ge \left(1/2+\alpha\right)n\ge \left(1/2+ \alpha/2\right)m$
for sufficiently large $n$.  
Hence, by~\cref{lem:normal-fpm}, with $\alpha/2$ in place
of $\alpha$, there exists a constant $C>1$ and a $C$-normal perfect fractional matching $\omega_{\mathbf s}:E(L_H(\mathbf s))\to\mathbb{R}_{\ge0}.$
We choose such a fractional matching for every ordered
$(k-2)$-tuple $\mathbf s$. 
When several $k$-graphs are involved, these choices are made separately for each graph.
For $k=2$, the tuple $\mathbf s$ is empty and $L_H(\mathbf s)=H$.

Let $\mathbf x=(x_1,s_1,\ldots,s_{k-2})\in\Omega_{k-1}(V(H)).$
A \emph{balanced extension} of $\mathbf x$ in $H$ is obtained by choosing $y\in N_H(\{x_1,s_1,\ldots,s_{k-2}\})$ with probability
\begin{equation}\label{eq:balanced extension}
  q_H(\mathbf x,y):=\omega_{\mathbf s}(x_1y),
\end{equation}
and replacing $\mathbf x$ by $\mathbf x\circ y=(s_1,\ldots,s_{k-2},y).$
This is indeed a probability distribution, since $\omega_{\mathbf s}$ is a perfect fractional matching, hence $\sum_{y\in N_{L_H(\mathbf s)}(x_1)}\omega_{\mathbf s}(x_1y)=1.$

We shall use the following two elementary properties of balanced
extensions.

\begin{lemma}\label{lem:property of banlanced extension}
    Let $1/n\ll\alpha,1/k$, where $k\ge2$. 
    Let $H$ be an $n$-vertex $k$-graph satisfying $\delta_{k-1}(H)\ge(1/2+\alpha)n$.
    Then there exists $c_0=c_0(k,\alpha), C_0=C_0(k,\alpha)>0$ such that the following hold. 
    \begin{enumerate}[$(i)$]
        \item  For every ordered $(k-1)$-tuple $\mathbf{x}$ and every legal balanced extension $y$ of $\mathbf x$, $\frac{c_0}{n}\le q_H(\mathbf x,y)\le \frac{C_0}{n}.$
        \item Suppose that $\mathbf X$ is uniformly distributed on $\Omega_{k-1}(V(H))$.
        Given $\mathbf X$, let $\mathbf X^+$ be the resulting ordered $(k-1)$-tuple after performing one balanced extension in $H$ conditioned on $\mathbf X$.
        Then $\mathbf X^+$ is also uniformly distributed on $\Omega_{k-1}(V(H))$.
    \end{enumerate}
\end{lemma}

\begin{proof}
    $(i)$ is immediate by the definition of balanced extension.
    For $(ii)$, fix a target state $(\mathbf{s},y)\in \Omega_{k-1}(V(H))$, where $\mathbf{s}\in \Omega_{k-2}(V(H))$.
    Then 
    \begin{align*}
        \Pr(\mathbf{X}^+=(\mathbf{s},y))
        &=\sum_{x\in N_{L_H(\mathbf s)}(y)}\Pr(\mathbf{X}^+=(\mathbf{s},y),\mathbf{X}=(x,\mathbf{s}))\\
        % &=\sum_{x\in N_{L_H(\mathbf s)}(y)}\Pr(\mathbf{X}^+=(\mathbf{s},y)\mid\mathbf{X}=(x,\mathbf{s}))\Pr(\mathbf{X}=(x,\mathbf{s}))\\
        &=\frac{1}{|\Omega_{k-1}(V(H))|}\sum_{x\in N_{L_H(\mathbf s)}(y)}\omega_{\mathbf s}(xy)=\frac{1}{|\Omega_{k-1}(V(H))|},
    \end{align*}
    where the last equality follows from the perfectness of $\omega_{\mathbf s}$ at the vertex $y$.
    Thus $\mathbf X^+$ is uniformly distributed on $\Omega_{k-1}(V(H))$.
\end{proof}

\subsection{Forgetting the initial state}

We next show that repeated balanced extensions rapidly erase the dependence on the initial ordered tuple.
This will allow us to reach any positive-density family of terminal tuples with uniformly positive probability.

For probability distributions $p$ and $q$ on a finite set $\Omega$,
a \emph{coupling} of $p$ and $q$ is a pair of random variables
$(X,Y)$ defined on the same probability space, with marginal
distributions $p$ and $q$, respectively. 
A coupling is called \emph{maximal} if it maximizes $\Pr(X=Y)$. 
It is a standard fact that a maximal coupling exists and satisfies 
$\Pr(X=Y)=\sum_{z\in\Omega}\min\{p(z),q(z)\}.$
See ~\cite[Section~4.2 and Chapter~5]{levin2026markov}.

\begin{lemma}\label{lem:initial state}
    Let $1/n\ll\alpha,1/k$, where $k\ge 2$.
    There exists $\rho=\rho(k,\alpha)>0$ such that the following holds. Let $H_1,H_2,\dots$ be a sequence of $k$-graphs on a vertex set $V=[n]$ satisfying $\delta_{k-1}(H_i)\ge (1/2+\alpha)n$.
    Start from an arbitrary $\mathbf{x}\in \Omega_{k-1}(V)$ and let $\mathbf{X}_0=\mathbf{x}$.
    For each step $i\ge 1$, let $\mathbf{X}_i$ be the tuple obtained by performing a balanced extension in $H_i$ conditioned on $\mathbf{X}_{i-1}$.
    Then for every $t\ge 0$ and every $\mathcal{B}\subseteq \Omega_{k-1}(V)$,
    \[
    \left|\Pr(\mathbf X_t\in \mathcal{B})-\frac{|\mathcal{B}|}{|\Omega_{k-1}(V)|} \right|\le (1-\rho)^{\lfloor\frac{t}{k-1}\rfloor}.
    \]
\end{lemma}
\begin{proof}
    We couple two balanced walks $(\mathbf X_i)_{i\ge 0}$ and $(\mathbf Y_i)_{i\ge 0}$.
    The first starts from a given state $\mathbf X_0=\mathbf x\in \Omega_{k-1}(V)$, whereas $\mathbf Y_0$ is uniformly distributed on $\Omega_{k-1}(V)$.
    At step $i$, both walks perform a balanced extension in $H_i$.
    By ~\cref{lem:property of banlanced extension}$(ii)$, $\mathbf Y_i$ is uniformly distributed on $\Omega_{k-1}(V)$ for every $i\ge0$.

    We now specify the coupling of the two transitions. 
    At some step, suppose that the current states are $\mathbf u=(u_1,\ldots,u_{k-1})$ and $\mathbf v=(v_1,\ldots,v_{k-1})$
    and that the next extension is performed in a $k$-graph $H$.
    Let $p(z):=q_H(\mathbf u,z)$ and $q(z):=q_H(\mathbf v,z)$ be the corresponding transition distributions. 
    Writing $U:=\{u_1,\ldots,u_{k-1}\}$, $W:=\{v_1,\ldots,v_{k-1}\}$, we have $|N_H(U)\cap N_H(W)|\ge d_H(U)+d_H(W)-n\ge 2\alpha n$.
    For every $z\in N_H(U)\cap N_H(W),$ ~\cref{lem:property of banlanced extension}$(i)$ gives $p(z),q(z)\ge \frac{c_0}{n}.$
    Consequently, $\sum_{z\in V}\min\{p(z),q(z)\}\ge|N_H(U)\cap N_H(W)|c_0/n\ge 2\alpha c_0.$
    We couple the two extensions maximally. 
    Thus, conditional on the current states, the two walks append the same vertex with probability at least $2\alpha c_0$.

    Now we partition the evolution into consecutive blocks of $k-1$ steps. Conditional on the history up to the beginning of any such block, the probability that the two walks append the same vertex at every step of the block is at least $\rho:=(2\alpha c_0)^{k-1}$.
    If this event occurs, then after these $k-1$ extensions the two ordered tuples coincide, since all entries of the initial states have been shifted out. 
    Once the two walks meet, we couple all subsequent extensions identically, so that they remain together forever.
    Thus, we have $ \Pr(\mathbf X_t\ne\mathbf Y_t)\le(1-\rho)^{\lfloor t/(k-1)\rfloor}.$
    Hence, for every $\mathcal B\subseteq\Omega_{k-1}(V)$,
    \[
    \left|\Pr(\mathbf X_t\in\mathcal B)-\frac{|\mathcal B|}{|\Omega_{k-1}(V)|}\right|
    =\left|\Pr(\mathbf X_t\in\mathcal B)-\Pr(\mathbf Y_t\in\mathcal B)
    \right|
    \le\Pr(\mathbf X_t\ne\mathbf Y_t)
    \le
    (1-\rho)^{\lfloor t/(k-1)\rfloor},
    \]
    as required.
\end{proof}

We are now ready to prove the connecting lemma. 
We run the balanced walk for sufficiently many steps and use ~\cref{lem:initial state} to ensure that its terminal state falls into a large family of tuples from which the prescribed terminal end can be reached.

\begin{proof}[Proof of ~\cref{lem:connecting}]
    Let $\mathbf{x}=(x_1,\dots,x_{k-1})$ and $\mathbf{y}=(y_1,\dots,y_{k-1})$ be the prescribed tuples.
    Since $|R\setminus R'|\le \alpha N/2$, for every $i\in[r+k-1]$ and every $(k-1)$-set $S$, $|N_{J_i}(S)\cap R'|\ge (1/2+\alpha/2)N$. 
    Set $U:=R'\cup V(\mathbf x)\cup V(\mathbf y).$
    Then $|R'|\le |U|\le |R'|+2(k-1)\le N+2(k-1).$
    In particular, for every $i\in[r]$ and every $S\in\binom{U}{k-1}$,
    $d_{J_i[U]}(S)\ge (1/2+\alpha/2)N \ge (1/2+\alpha/3)|U|$ for sufficiently large $N$.
    
    Let $\mathcal{B}$ be a set of ordered $(k-1)$-tuples $\mathbf{z}=(z_1,\dots,z_{k-1})\in\Omega_{k-1}(R')$, which are disjoint from $V(\mathbf{x})\cup V(\mathbf{y})$ and for which $z_1\dots z_{k-1}y_1\dots y_{k-1}$ is a tight path using the colors $J_{r+1},\dots,J_{r+k-1}$.

    We construct $\mathbf{z}$ backward.
    First choose $z_{k-1}\in R'$ such that $\{z_{k-1},y_1,\dots,y_{k-1}\}\in E(J_{r+k-1})$, next choose $z_{k-2}\in R'$ such that $\{z_{k-2},z_{k-1},y_1,\dots,y_{k-2}\}\in E(J_{r+k-2})$ and continue in this way.
    At each step, there are at least $(1/2+\alpha/2)N-O_k(1)\ge N/2$ available choices after excluding the vertices already chosen and the vertices of $V(\mathbf x)\cup V(\mathbf y)$.
    Consequently, $|\mathcal{B}|\ge (N/2)^{k-1}$. 
    Since $|U|\le N+2(k-1)$, there is a constant $c_1=c_1(k)>0$ such that $|\mathcal{B}|/|\Omega_{k-1}(U)|\ge c_1(k)>0$.
    
    Starting from $\mathbf{x}$, perform a balanced extension in $J_i[U]$ at step $i$, for $i=1,\dots,r$, and let $\mathbf X_t$ denote the state after $t$ extensions.
    Let $\rho=\rho(k,\alpha/3)>0$ be the constant supplied by \cref{lem:initial state}.
    Since $1/r\ll\alpha,1/k$, we may assume that $(1-\rho)^{\lfloor r/(k-1)\rfloor}\le c_1/4$.
    By ~\cref{lem:initial state}, $\Pr(\mathbf X_r\in\mathcal{B})\ge c_1 -(1-\rho)^{\lfloor r/(k-1)\rfloor}\ge 3c_1/4$.

    It remains to ensure that the $r$ vertices appended during the walk are distinct belong to $R'$, and avoid the two prescribed ends. 
    At every step there are at most $2(k-1)+r$ forbidden vertices, namely the vertices of $V(\mathbf x)\cup V(\mathbf y)$ together with those already appended.
    By ~\cref{lem:property of banlanced extension}$(i)$, applied to the
    graphs $J_i[U]$, the conditional probability of choosing a forbidden
    vertex at any given step is $O_{k,\alpha,r}(1/N).$ 
    Since $r$ is fixed, a union bound over all $r$ steps shows that the probability that a forbidden vertex is chosen is at most $O_{k,\alpha,r}(1/N)\le c_1/4$.
    Hence, with probability at least $c_1/2$, we obtain a desired $\mathcal{J}$-colored tight path with ends $\mathbf{x}$ and $\mathbf{y}$.
    
    Finally, we have every fixed ordered $r$-tuple $(\omega_1,\dots,\omega_r)$ is generated with probability at most $\left(C_0/|U|\right)^r\le\left(C_0/|R'|\right)^r$, where $C_0=C_0(k,\alpha)$ is a constant.
    Let $\mathcal{S}$ be the family of $r$-sets of $R'$ which occur as the internal vertex set of such an exactly $\mathcal J$-colored $\mathbf x,\mathbf y$-path.
    Each member of $\mathcal S$ has at most $r!$ orderings, 
    hence $ c_1/2\le r!|\mathcal S|\left(C_0/|R'|\right)^r.$
    Therefore $|\mathcal S|\ge\frac{c_1}{2r!C_0^r}|R'|^r=\Omega(|R'|^r),$ as required.
\end{proof}

\section{Swap absorption}\label{sec:absorption}

In this section, we prove ~\cref{lem:absorber}. 
Our absorption argument is inspired by the swap-absorption method for Hamilton cycles originating from Reiher, R{\"o}dl, Ruci{\'n}ski, Schacht and Szemer{\'e}di~\cite{reiher2019minimum}, and by the reachability approach of Lo and Markstr{\"o}m ~\cite{lo2015f}, originally in the context of hypergraph matchings and tilings.

Our absorbing argument must preserve not only the ordered ends of a path but also its prescribed sequence of colors.
We therefore work with path blocks, each of which specifies two ordered ends and a fixed color list.
The associated auxiliary hypergraph saves the internal vertex sets that realize this prescribed path.
The key observation is that any two vertices admit many common completions in this auxiliary hypergraph. 
Consequently, one vertex can be exchanged for another inside a block without changing its ends or color list.
We combine these local exchanges through a collection of randomly selected blocks, obtaining an absorber that works simultaneously for every prescribed subset of the reservoir of the required size.
First, we provide the necessary definitions.

\begin{definition}[Path block]
    Let $\mathcal J=(J_1,\ldots,J_{q+k-1})$ be an ordered collection of $k$-graphs on a common vertex set $V$.
    A \emph{$q$-block} $\mathcal B=(\mathbf x,\mathbf y;\mathcal J)$ consists of two disjoint ordered $(k-1)$-tuples $\mathbf x,\mathbf y\in \Omega_{k-1}(V)$, called the \emph{ends} of
    $\mathcal B$, together with the prescribed color list
    $\mathcal J$.
    The \emph{auxiliary $q$-graph} $F_{\mathcal B}$ has vertex set
    $ V(F_{\mathcal B}):=V\setminus\bigl(V(\mathbf x)\cup V(\mathbf y)\bigr).$
    A set $Q\in\binom{V(F_{\mathcal B})}{q}$ is an edge of $F_{\mathcal B}$ if there exists an ordering $Q=(z_1,\ldots,z_q)$ such that $x_1\cdots x_{k-1}z_1\cdots z_qy_1\cdots y_{k-1}$ is an exactly $\mathcal J$-colored tight $\mathbf x,\mathbf y$-path.
    In this case, we call $Q$ a solution of $F_{\mathcal{B}}$.
\end{definition}

Thus, the edges of $F_{\mathcal B}$ are precisely the possible internal vertex sets of an exactly $\mathcal J$-colored tight path with prescribed ends $\mathbf x,\mathbf y$.
For $k=2$, a $q$-block consists of two prescribed vertices and a list of $q+1$ colors.

The block formulation is useful beyond the proof of the absorbing-path lemma.
In particular, it will be used again in ~\cref{subsec:tree with few leaves}, where absorbing blocks are distributed among several prescribed bare paths.

\begin{definition}[$c$-reachable]
Let $\mathcal B$ be a $q$-block, let $R\subseteq V(F_{\mathcal B})$ with $|R|=N$, and let $c>0$. 
Two distinct vertices $u,v\in V(F_{\mathcal B})$ are \emph{$c$-reachable in $\mathcal B$ with respect to $R$} if there are at least $cN^{q-1}$ sets $S\in\binom{R\setminus\{u,v\}}{q-1}$ such that
$S\cup\{u\},\ S\cup\{v\}\in E(F_{\mathcal B}).$
Each such $S$ is called a \emph{reachable set} for $u$ and $v$ in $\mathcal B$.
\end{definition}

In particular, if $S$ is a reachable set for $u$ and $v$ in $\mathcal{B}$, then $\mathcal B$ admits a solution with internal vertex set $S\cup\{u\}$ and another with internal vertex set $S\cup\{v\}$.
The two corresponding paths need not induce the same ordering on $S$; what matters is that they have the same ordered ends and use the same prescribed color list.

Our next lemma shows that every block has many solutions and that every pair of vertices is $c$-reachable under the codegree condition.

\begin{lemma}\label{lem:block-reachability}
Let $0<1/N\ll\alpha,1/k$, where $k\ge 2$. 
There exist $q=q(k,\alpha)\in\mathbb N$ and $c=c(k,\alpha)>0$ such that the following holds. 
Let $\mathcal B=(\mathbf x,\mathbf y;\mathcal J)$ be a $q$-block, where $\mathcal J=(J_1,\ldots,J_{q+k-1})$, and let $R\subseteq V(F_{\mathcal B})$ be a $N$-set satisfying $\delta_{k-1}(\mathcal{J},R)\ge (1/2+\alpha)N$. 
Then, every two distinct vertices $a,b\in V(F_{\mathcal B})$ are $c$-reachable in $\mathcal B$ with respect to $R$, regardless of whether
$a$ and $b$ belong to $R$. 
Moreover, $|E(F_{\mathcal B}[R])|\ge cN^q.$
\end{lemma}
\begin{proof}
    Choose $r=r(k,\alpha)$ sufficiently large from ~\cref{lem:connecting}, and set $q:=2r+2k-1$. 
    We first prove the reachability assertion. 
    Fix distinct $a,b\in V(F_{\mathcal B})$. 
    We shall construct a central segment $u_1\cdots u_{k-1}\,z\,w_1\cdots w_{k-1}$ which uses the prescribed $k$ middle colors for both choices $z=a$ and $z=b$.

    Choose distinct vertices $u_1,\ldots,u_{k-2}\in R\setminus(\{a,b\}\cup V(\mathbf x)\cup V(\mathbf y)).$
    There are $\Omega(N^{k-2})$ ways to do this.
    For $\xi\in\{a,b\}$, put $A_\xi:=\{u_1,\ldots,u_{k-2},\xi\}.$
    Since $|N_{J_{r+k}}(A_a)\cap R|,|N_{J_{r+k}}(A_b)\cap R|\ge(1/2+\alpha)N,$ we have $|N_{J_{r+k}}(A_a)\cap N_{J_{r+k}}(A_b)\cap R|\ge 2\alpha N.$
    Hence there are at least $2\alpha N-O_k(1)$ choices for a vertex $u_{k-1}$ such that $u_1\cdots u_{k-1}\xi$ forms an edge of $J_{r+k}$ for both $\xi=a$ and $\xi=b$.
    We now choose $w_1,\ldots,w_{k-1}$ successively. 
    More precisely, having chosen $w_1,\ldots,w_{j-1}$, where $1\le j\le k-1$, consider for $\xi\in\{a,b\}$ the $(k-1)$-set
    $A_\xi^{(j)}:=\{u_{j+1},\ldots,u_{k-1},\xi,w_1,\ldots,w_{j-1}\},$
    where the first part is empty when $j=k-1$.
    By the minimum codegree assumption, $|N_{J_{r+k+j}}(A_a^{(j)})\cap R|, |N_{J_{r+k+j}}(A_b^{(j)})\cap R|\ge(1/2+\alpha)N.$
    Thus their intersection contains at least $2\alpha N$ vertices.
    After excluding the bounded number of vertices already used, there are still $\Omega(N)$ choices for $w_j$.
    It follows that there are $\Omega(N^{2k-2})$ ordered choices of $\mathbf u=(u_1,\ldots,u_{k-1}),$ $\mathbf w=(w_1,\ldots,w_{k-1})$ 
    such that, for each $\xi\in\{a,b\}$, $u_1\cdots u_{k-1}\xi  w_1\cdots w_{k-1}$ is a tight path whose $k$ consecutive edges use the colors $J_{r+k},J_{r+k+1},\ldots,J_{r+2k-1}$ in order.

    Fix one such choice of $\mathbf u$ and $\mathbf w$. 
    We next join $\mathbf x$ to $\mathbf u$ using the first $J_1,\ldots,J_{r+k-1}$ and join $\mathbf w$ to $\mathbf y$ using the last $J_{r+2k},\ldots,J_{q+k-1}.$
    For the left connector, apply \cref{lem:connecting} with its internal vertices chosen in $R\setminus
    \bigl(V(\mathbf u)\cup V(\mathbf w)\cup\{a,b\}\cup V(\mathbf x)\cup V(\mathbf y)\bigr).$
    Only $O_k(1)$ vertices have been excluded, so we have $\Omega(N^r)$ choices.
    For each such choice, we apply \cref{lem:connecting} once
    more to obtain the right connector, additionally avoiding the $r$
    internal vertices of the left connector. 
    Again only $O_{k,r}(1)$ vertices of $R$ are forbidden, and hence there are $\Omega(N^r)$ choices for the right connector.

    Let $S$ consist of all internal vertices used in these two connectors and in $\mathbf u,\mathbf w$.
    Thus, $|S|=r+(k-1)+(k-1)+r=2r+2k-2=q-1,$ and $S\subseteq R\setminus\{a,b\}.$
    Hence, both $S\cup\{a\}$ and $S\cup\{b\}$ are solutions of the same block $\mathcal B$. 
    Hence $S$ is a reachable set for $a$ and $b$.
    The number of ordered choices is $\Omega(N^{2k-2})\cdot
    \Omega(N^r)\cdot\Omega(N^r)=\Omega(N^{2r+2k-2})=\Omega(N^{q-1})$.
    A fixed $(q-1)$-set $S$ is counted at most $(q-1)!$ times. 
    Therefore there is a constant $c_1=c_1(k,\alpha)>0$ such that at least
    $c_1N^{q-1}$ distinct reachable sets exist for every pair of distinct vertices $a,b\in V(F_{\mathcal B})$.
    
    To obtain the density assertion, fix $a\in R$ and choose $b\in R\setminus\{a\}$.
    Every reachable set $S$ for $a,b$ yields a distinct edge
    $S\cup\{a\}\in E(F_{\mathcal B}[R])$, so $d_{F_{\mathcal B}[R]}(a)\ge c_1N^{q-1}.$
    Summing over $a\in R$ and counting each edge exactly $q$ times gives $q|E(F_{\mathcal B}[R])|=\sum_{a\in R}d_{F_{\mathcal B}[R]}(a)\ge c_1N^q.$
    Taking $c:=c_1/q$ proves both assertions.
\end{proof}

We now assemble these local swaps into a global absorber. 
We first do so for an arbitrary family of blocks, and the concatenation will give ~\cref{lem:absorber}. 
Keeping the blocks separate will also be useful for trees.

\begin{lemma}\label{lem:block absorb}
Let $0<1/n\ll\beta\ll\theta\ll\alpha,1/k$, where $k\ge2$. 
There is an integer $q=q(k,\alpha)$ such that the following holds. 
Put $b=\lfloor\beta n\rfloor$ and $T=\lceil b/\theta\rceil+b$. 
Let $\mathcal{J}$ be a collection of $k$-graphs on $V=[n]$ with $\delta_{k-1}(\mathcal{J})\ge(1/2+\alpha)n$, and let $Z\subseteq V$ satisfy $b\le|Z|\le2b$. 
Given $T$ prescribed $q$-blocks with color lists in $\mathcal{J}$ and ends outside $Z$, let $D$ be the union of their ends. 
Then there is $A_0\subseteq V\setminus(Z\cup D)$, of size $qT-b$ such that for every $Z'\in\binom Zb$, the blocks have pairwise disjoint solutions whose union is $A_0\cup Z'$.
\end{lemma}
Here the ends of different blocks are not required to be disjoint.
Only their internal vertex sets will be chosen pairwise disjoint.
\begin{proof}
    Take $q,c$ from ~\cref{lem:block-reachability} with $\alpha/2$ in place of $\alpha$ and reduce $c$ if necessary such that $c\le1$.
    Choose constants $1/n\ll\beta\ll\theta\ll p\ll \alpha,c,1/q,1/k$.
    Put $M=\lceil b/\theta\rceil$, so that $T=M+b$.
 Let $F_t$ be the auxiliary $q$-graph of block $t$, for every $t\in[T]$. 
 Call blocks $1,\ldots,M$ normal and blocks $M+1,\ldots,M+b$ absorbed. 
 We have $|D|\le2T(k-1)$.

For $i=1,\ldots,b$,  we will choose a fixed $(q-1)$-set $C_i$ outside $Z\cup D\cup\bigcup_{j<i}C_j$. 
At step $i$, let $R_i=V\setminus(Z\cup D\cup\bigcup_{j<i}C_j)$ and $N_i:=|R_i|$. 
By $\beta\ll\theta\ll\alpha,c$, we have $|Z|+|D|+(q-1)b+qM\le\min\{\alpha/10,c/4\}n$.
In particular, $N_i\ge n/2$.
Moreover, deleting $V\setminus R_i$ reduces every codegree by at most
$\alpha n/10$, and hence $\delta_{k-1}(\mathcal J,R_i)\ge (1/2+\alpha/2)N_i.$
By~\cref{lem:block-reachability}, $|E(F_{M+i}[R_i])|\ge cN_i^q.$
For each $C\in\binom{R_i}{q-1}$, let $L(C):=\{v\in R_i\setminus C:C\cup\{v\}\in E(F_{M+i})\}$.
Since $\sum_{C\in\binom{R_i}{q-1}}|L(C)|=q|E(F_{M+i}[R_i])|\ge qcN_i^q,$ the pigeonhole principle yields a $(q-1)$-set
$C_i\subseteq R_i$ such that $|L(C_i)|\ge cN_i\ge cn/2.$
Set $L_i:=L(C_i)$.
After choosing all the sets $C_i$, the total number
of vertices removed by these sets and the subsequent
normal configurations is at most $(q-1)b+qM\le cn/4.$
Consequently, every subsequent available set contains at least $cn/4$ vertices of each $L_i$.

We next choose the normal configurations in order. 
Suppose that configurations $W_1,\dots,W_{g-1}$ have already been chosen, and put $R_g=V\setminus(Z\cup D\cup\bigcup_i C_i\cup\bigcup_{h<g}W_h)$. 
We sample a uniformly random injection $\varphi:[q]\to R_g$, conditioned on $\varphi([q])\in E(F_g)$.
Put $v_g=\varphi(1)$, $S_g=\{\varphi(2),\ldots,\varphi(q)\}$ and $W_g=S_g\cup\{v_g\}$. 
This conditional distribution is well-defined, since ~\cref{lem:block-reachability} guarantees that $|E(F_g[R_g])|\ge c|R_g|^q$.
For $i\in[b]$ and $u\in Z$, we call the $g$th block an $(i,u)$-swap absorber if $v_g\in L_i$ and $S_g\cup\{u\}\in E(F_g)$.
We claim that for every $i\in[b]$ and $u\in Z$, the conditional probability that the $g$th block is an $(i,u)$-swap absorber is at least $p$, regardless of the previous choices.
To prove this, fix $i\in[b]$ and $u\in Z$, and condition on an arbitrary history $\mathcal F_{g-1}$ of the preceding choices.
Let $\mathcal G_g$ be the event $W_g\in E(F_g)$ and $\mathcal E_g(i,u)$ be the event that the unconditioned trial produces an $(i,u)$-swap absorber.

There are at least $cn/4$ choices for $v_g\in L_i\cap R_g$.
For each $(v_g,u)$, ~\cref{lem:block-reachability} gives at least $c|R_g|^{q-1}$ reachable sets $S_g$. 
Hence, for a uniformly random injection before conditioning on $\mathcal G_g$, we have $\mathbb \Pr(\mathcal E_g(i,u)\mid\mathcal F_{g-1})\ge(cn/4)c|R_g|^{q-1}/|R_g|^q\ge c^2/4\ge p,$
Since $\mathcal E_g(i,u)\subseteq\mathcal G_g$, we have 
$\Pr(\mathcal E_g(i,u)\mid
\mathcal F_{g-1},\mathcal G_g)=\Pr(\mathcal E_g(i,u)\mid\mathcal F_{g-1})/\Pr(\mathcal G_g\mid\mathcal F_{g-1})\ge p.$
Thus, regardless of the previous choices, the conditional probability that the $g$th normal block is an $(i,u)$-swap absorber is at least $p$.

\begin{claim}\label{claim:swaps}
With positive probability, every target $(i,u)\in[b]\times Z$ has at least $pM/2$ $(i,u)$-swap absorbers in normal blocks.
\end{claim}
\begin{proof}
Fix $(i,u)$ and let $X_g$ be the indicator of the event that the
$g$th normal block is an $(i,u)$-swap absorber. 
Conditional on the previous history $\mathcal{F}_{g-1}$, although the random variables $X_1,\ldots,X_M$ need not be independent, 
    \[
    \mathbb E(2^{-X_g}\mid\mathcal F_{g-1})=\Pr(X_g=0\mid\mathcal F_{g-1})+\frac12\Pr(X_g=1\mid\mathcal F_{g-1})\le 1-\frac p2.
    \]
    Iterating conditional expectations gives 
    $\mathbb{E}\left(2^{-\sum_{g=1}^{M}X_g}\right)
    \le \left(1-p/2\right)^M\le e^{-pM/2}.$
    Hence, by Markov's inequality,
    \[
        \Pr\left(\sum_{g=1}^{M}X_g<\frac{pM}{2}\right)=\Pr\left(2^{-\sum_{g=1}^{M}X_g}>2^{-\frac{pM}{2}}\right)
        \le 2^{pM/2}e^{-pM/2}
        \le e^{-pM/8}.
    \]
    
    There are at most $n^2$ choices of $(i,u)$, while $M=\Theta(n)$.
     A union bound therefore shows that, with positive probability, the desired property holds simultaneously for every target $(i,u)$.
\end{proof}
We fix such a realization and put $A_0=\bigcup_i C_i\cup\bigcup_g W_g$. 
These sets are disjoint and $|A_0|=(q-1)b+qM=qT-b$. 
Now let $Z'=\{u_1,\ldots,u_b\}\subseteq Z$.
Since $b\le\theta M<pM/2$, we can choose distinct normal blocks $g_1,\ldots,g_b$ such that $g_i$ is an $(i,u_i)$-swap absorber. 
In block $g_i$, we replace $S_{g_i}\cup\{v_{g_i}\}$ by $S_{g_i}\cup\{u_i\}$ and fill block $M+i$ with $C_i\cup\{v_{g_i}\}$. 
The solutions in all other normal blocks remain unchanged.
These solutions are pairwise disjoint and use exactly $A_0\cup Z'$, without changing any ends or color lists. 
This proves the lemma.
\end{proof}

\begin{proof}[Proof of ~\cref{lem:absorber}]
    Choose $1/n\ll\beta\ll\theta\ll\alpha,1/k$, and let $q$ be given by ~\cref{lem:block absorb}.
    Put $b=\lfloor\beta n\rfloor, T=\lceil b/\theta\rceil+b$ and $a=T(q+k-1)\le\alpha n/20$.
    Choose pairwise disjoint ordered $(k-1)$-tuples $D_0,\ldots,D_T$ outside $Z$. 
    For each $t\in[T]$, let block $t$ have ends $D_{t-1},D_t$ and color list $J_{(t-1)(q+k-1)+1},\ldots,J_{t(q+k-1)}$. 
    Apply ~\cref{lem:block absorb} to obtain $A_0$. 
    Set $A=A_0\cup\bigcup_{t=0}^T V(D_t)$, $\mathbf{x}=D_0$ and $\mathbf{y}=D_T$. 
    Then $|A|=qT-b+(T+1)(k-1)=a+k-1-b$.
    For every $Z'\in\binom Zb$, concatenate those blocks in their block order. 
    They form an exactly $\mathcal{J}$-colored tight $\mathbf x,\mathbf y$-path on $A\cup Z'$, as required.
\end{proof}

\section{Almost path cover}\label{sec:almost path cover}

In this section, we prove ~\cref{lem:almost cover}. 
We use a random greedy process which selects the required tight paths one at a time while maintaining the codegree condition in the unused vertex set.

\begin{proof}[Proof of ~\cref{lem:almost cover}]
    For $t\in\{0,\ldots,s\}$, put $m_t:=n-Kt.$ 
    Whenever the process reaches time $t$, let $R_t$ denote the set of unused vertices after $P_1,\ldots,P_t$ have been chosen. 
    Thus $|R_t|=m_t$.
    Since $s=\left\lfloor(1-\varepsilon)n/K\right\rfloor,$ we have $m_t\ge \varepsilon n$ for every $t\le s.$

    For later use, let $\mathcal N:=\left\{N_{J_c}(S):c\in[m],S\in\binom{V}{k-1}\right\}.$
    We continue the process as long as 
    \begin{equation}\label{eq:almost-cover-good}
    |T\cap R_t|
    \ge (1/2+\alpha/2)m_t
    \qquad\text{for every }T\in\mathcal N.
    \end{equation}
    In particular, for every $c\in[m]$ and every $S\in\binom{R_t}{k-1}$,
    $d_{J_c[R_t]}(S)=|N_{J_c}(S)\cap R_t|\ge (1/2+\alpha/2)m_t.$
    Thus balanced extensions are available in every graph that may occur in the prescribed color pattern.
    
    Suppose that \eqref{eq:almost-cover-good} holds at time $t<s$. 
    We now describe the choice of $P_{t+1}$. 
    We generate a trial by first choosing $(x_1,\ldots,x_{k-1})$ uniformly from $\Omega_{k-1}(R_t)$ and then performing balanced extensions in $J_{\psi_{t+1}(j)}[R_t]$ for $j=1,\ldots,K-k+1$.
    This produces a sequence $x_1,\ldots,x_K$, which may contain repeated vertices.

    Before conditioning on distinctness, every one of the $K$ positions is uniformly distributed on $R_t$. 
    Indeed, this is clear for the initial $(k-1)$-tuple, and ~\cref{lem:property of banlanced extension}$(ii)$ shows inductively that after each extension the current ordered $(k-1)$-tuple remains uniform on $\Omega_{k-1}(R_t)$.

    Let $\mathcal E_t$ be the event that $x_1,\ldots,x_K$ are pairwise distinct.
    By ~\cref{lem:property of banlanced extension}$(i)$, applied with $\alpha/2$ in place of $\alpha$, the conditional probability that an appended vertex coincides with any previously exposed vertex is $O_{k,K,\alpha}(1/m_t)$. 
    Hence $\Pr(\mathcal E_t^c\mid\mathcal F_t)
    =O_{k,K,\alpha}(1/m_t)$, where $\mathcal F_t$ denotes the history up to time $t$.

    We repeatedly generate an independent trial until $\mathcal E_t$ occurs, and then let $P_{t+1}:=x_1\cdots x_K.$
    Thus $P_{t+1}$ is a $K$-vertex tight path which is exactly $\psi_{t+1}$-colored. 
    We then set $R_{t+1}:=R_t\setminus V(P_{t+1}).$

    We next show that the process does not stop before $s$ with positive probability. 
    Fix $T\in\mathcal N$, and write $d_t(T):=|T\cap R_t|$, $Y_t(T):=d_t(T)/m_t.$
    Initially, $Y_0(T)\ge 1/2+\alpha.$
    Let $X_t(T):=|T\cap V(P_{t+1})|.$
    We claim that, whenever the process has not yet stopped, we have \begin{equation}\label{eq:almost-cover-expectation}
    \mathbb E\bigl[X_t(T)\mid\mathcal F_t\bigr]=K Y_t(T)+O_{k,K,\alpha}(m_t^{-1}).
    \end{equation}

    For an unconditioned trial, define $L_t(T):=\sum_{j=1}^{K}\mathbf 1_{\{x_j\in T\}}.$
    Since each of the $K$ positions is uniform on $R_t$, $\mathbb E\bigl[L_t(T)\mid\mathcal F_t\bigr]=Kd_t(T)/m_t=KY_t(T).$
    On the event $\mathcal E_t$, the vertices $x_1,\ldots,x_K$ are pairwise distinct, so $L_t(T)=|T\cap\{x_1,\ldots,x_K\}|.$
    Moreover, since the accepted trial has the same distribution as an unconditioned trial conditioned on $\mathcal E_t$, we obtain
    $\mathbb E\bigl[X_t(T)\mid\mathcal F_t\bigr]=
    \mathbb E\bigl[L_t(T)\mid\mathcal F_t,\mathcal E_t\bigr].$
    Since $0\le L_t(T)\le K$, conditioning on $\mathcal E_t$ changes its expectation by at most
    $\left|\mathbb E[L_t(T)\mid\mathcal F_t,\mathcal E_t]-\mathbb E[L_t(T)\mid\mathcal F_t]\right|\le K\Pr(\mathcal E_t^c\mid\mathcal F_t)=O_{k,K,\alpha}(m_t^{-1})$.
    This proves \eqref{eq:almost-cover-expectation}.

    Since $m_{t+1}=m_t-K$ and $d_{t+1}(T)=d_t(T)-X_t(T),$ we have $Y_{t+1}(T)-Y_t(T)=\frac{m_tY_t(T)-X_t(T)}{m_{t+1}}-Y_t(T)=\frac{KY_t(T)-X_t(T)}{m_{t+1}}.$
    By~\eqref{eq:almost-cover-expectation}, we have $\left|\mathbb E\bigl[Y_{t+1}(T)-Y_t(T)\mid\mathcal F_t\bigr]\right|\le\frac{C}{m_tm_{t+1}}$ for some constant $C=C(k,K,\alpha)$.
    Moreover, since $0\le X_t(T)\le K$ and $0\le Y_t(T)\le1$, we have $|Y_{t+1}(T)-Y_t(T)|\le \frac{K}{m_{t+1}}.$

    Let $\tau$ be the first $t\le s$ for which \eqref{eq:almost-cover-good} fails, with $\tau=s+1$ if this never occurs. 
    For a fixed $T\in\mathcal N$, define $\widetilde Y_t(T):=Y_{\min\{t,\tau\}}(T)$.
    Thus, for every $0\le t<s$, $\widetilde Y_{t+1}(T)-\widetilde Y_t(T)=\mathbf 1_{\{t<\tau\}}\bigl(Y_{t+1}(T)-Y_t(T)\bigr).$
    Hence, by the estimates above, $\left|\mathbb E\left[\widetilde Y_{t+1}(T)-\widetilde Y_t(T)\mid\mathcal F_t\right]\right|\le\frac{C}{m_tm_{t+1}},$ and $\left| \widetilde Y_{t+1}(T)-\widetilde Y_t(T)\right|\le\frac{K}{m_{t+1}}.$
    Define
    \[
    M_t(T):=\widetilde Y_t(T)-Y_0(T)-\sum_{i=0}^{t-1}\mathbb E\left[\widetilde Y_{i+1}(T)-\widetilde Y_i(T)\mid\mathcal F_i\right].
    \]
    Then $(M_t(T))_{t=0}^{s}$ is a martingale.
    Moreover, $|M_{t+1}(T)-M_t(T)|\le\frac{K}{m_{t+1}}+\frac{C}{m_tm_{t+1}}\le\frac{C'}{m_{t+1}}$ for some constant $C'=C'(k,K,\alpha)$.
    Since $m_t\ge\varepsilon n$, $\sum_{t=0}^{s-1}\left(C'/m_{t+1}\right)^2=O(1/n).$
    Also, we have
    $\sum_{t=0}^{s-1}\left|\mathbb E\left[\widetilde Y_{t+1}(T)-\widetilde Y_t(T)\mid\mathcal F_t\right]\right|=O(1/n).$

    Suppose that $\tau\le s$, and that $T\in\mathcal N$ witnesses the failure at time $\tau$.
    Then $\widetilde Y_s(T)=Y_\tau(T)<1/2+\alpha/2,$ whereas $Y_0(T)\ge1/2+\alpha.$
    This implies, for sufficiently large $n$, $M_s(T)\le-\alpha/3.$
    By ~\cref{lem:Azuma}, $\Pr\bigl(M_s(T)\le-\alpha/3\bigr)\le\exp(-\Omega(n)).$
    Finally, $|\mathcal N|\le m\binom{n}{k-1}\le n^k.$
    A union bound therefore gives $\Pr(\tau\le s)\le n^k\exp(-\Omega(n))=o(1).$    
    Hence there is a realization for which the process continues through all $s$ steps. 
    In this realization $ P_1,\ldots,P_s$ are pairwise vertex-disjoint $K$-vertex tight paths, and $P_j$ is exactly $\psi_j$-colored for every $j\in[s]$, as required.
\end{proof}

\section{Spanning trees}\label{sec:tree}

In this section, we prove ~\cref{lem:leaves} and~\cref{lem:linkage}. 
Both proofs use the graph case of the tools established in ~\cref{sec:connecting,sec:absorption,sec:almost path cover}. 

\subsection{Trees with many leaves}
Suppose that the tree has linearly many leaves. 
We first remove a suitable collection of leaves and embed the remaining tree by a randomized procedure. 
During the embedding process, we will break down the tree into individual pieces and embed them one by one.
The deleted leaves will be embedded at the end using a perfect matching. 
The random embedding of each connected piece can be viewed as a tree-indexed random walk, an approach previously used for tree embedding by Sudakov and Vondr\'ak~\cite{sudakov2010randomized}.
The embedding must preserve every color neighborhood in the unused set and distribute the parents of the deleted leaves. 
These are the two degree conditions needed for the final perfect matching. 

\begin{proof}[Proof of ~\cref{lem:leaves}]
It suffices to prove the lemma for $0<\alpha<1/4$. 
For a fixed  pattern $\chi$, we discard all unused colors, so that at most $n-1$ color graphs remain.
Set $\eta=\alpha\lambda/1000$. 
Let $\rho\in(0,1)$ be given by ~\cref{lem:initial state} with $k=2$ and $\alpha/2$ in place of $\alpha$.
We choose a fixed integer $L\ge1$ with $(1-\rho)^L\le\eta$ and put $c=1/(4(L+2))$. 
For sufficiently large $n$, we write $D=\max\{2,\Delta(T)\}\le n^c$ and define $H=\lceil n^{1/4}\rceil,\ B=D^L,\ A=DH.$
Since $c\le1/12$, we have $A\le2n^{1/4+c}\le2n^{1/3},\ D^2A\le2n^{1/4+3c}\le2n^{1/2}.$

Choose $\ell=\lfloor\lambda n/2\rfloor\ge\lambda n/3$ leaves of $T$, denote their set by $\mathcal{L}$, and put $F=T-\mathcal{L}$. 
Fix a root of $F$. 
We show that we can partition $V(F)$ into connected pieces $F_1,\ldots,F_s$ such that each has order at most $A=DH$, with all but the root piece having at least $H$ vertices.

We obtain this partition from the leaves upwards. 
As long as at least $H$ vertices remain, we choose a deepest vertex whose remaining descendant subtree has at least $H$ vertices and remove that subtree as the next piece. 
Each child subtree has fewer than $H$ vertices, so the removed piece has at most $1+D(H-1)\le DH$ vertices. 
If the selected subtree is the whole remaining tree, stop; otherwise a final remainder of fewer than $H$ vertices forms the root piece. 
We then index the pieces from the root outwards. 
Every piece other than $F_1$ has exactly one edge to earlier pieces incident with its root, and $s\le n/H+1\le n^{3/4}+1$.

Call a vertex of a piece \emph{special} if its distance from that piece's root is less than $L$. 
There are at most $1+D+\cdots+D^{L-1}\le B$ special vertices per piece. 
The choice of $c$ gives $cL=1/4-2c$. 
Thus, we have $Bs<DBs\le2n^{1-c}=o(n)$.
We may assume $DBs\le\alpha\ell/12$.

We construct an embedding of $F$ into $V$ piece by piece, using $f$
to denote the partial embedding at each stage. 
The map $f$ is injective and satisfies $f(x)f(y)\in E(G_{\chi(xy)})$ for every edge $xy\in E(F)$ whose endpoints have been embedded. 
After embedding $F_1,\ldots,F_t$, let $R_t:=V\setminus f\left(\bigcup_{i=1}^{t}V(F_i)\right)$ be the set of unused host vertices. 
Write $a_t=|V(F_{t+1})|$, where $0\le t<s$ and let $m_t=n-\sum_{i=0}^{t-1}a_i.$
We have $a_t\le A$ and $m_t\ge m_s=\ell$. 
Let $\mathcal N=\{N_{G_i}(v):i\in\chi(E(T)),\ v\in V\}$, then $|\mathcal N|\le n(n-1).$
We continue as long as
\begin{equation}\label{eq:tree-invariant}
|W\cap R_t|\ge(1/2+\alpha/2)m_t
\qquad\text{for every }W\in\mathcal N.
\end{equation}
Here $\mathcal N$ includes neighborhoods of already used host vertices, since each new piece must attach to an embedded parent.  
Let $\mathcal F_t$ be the process through time $t$, keeping the history fixed after stopping.

At an active step $t<s$, the degree invariant allows us
to fix balanced-extension rules for every $G_i[R_t]$, using~\cref{lem:property of banlanced extension} with $\alpha/2$ in place of $\alpha$. 
For $t=0$, we choose the root image uniformly in $R_0$. 
For $t>0$, we choose it uniformly from the prescribed-color neighborhood of its already embedded parent inside $R_t$. 
This set is nonempty by \eqref{eq:tree-invariant}. 
We then expose the remaining vertices of $F_{t+1}$ in parent--child order, sampling each child's image according to the balanced-extension rule for its edge color.
Here, all vertices of $R_t$ are not removed throughout the trial.

Let $E_t$ be the event that all images in the trial are distinct. 
By ~\cref{lem:property of banlanced extension}(i), at each non-root exposure, any specified image has conditional probability at most $C/m_t$ for a constant $C=C(\alpha)$. 
A union bound over pairs of positions therefore gives
$\Pr(E_t^c\mid\mathcal F_t)\le C a_t^2/m_t\le 3C A^2/(\lambda n)=O(n^{-1/3}).$
We generate independent trials until $E_t$ occurs. 
The accepted trial gives an embedding of $F_{t+1}$, whose images we then remove from $R_t$.
The acceptance probability is positive, so this rejection step ends almost surely. 

Let $u\in V(F_{t+1})$ be non-special. 
The root-to-$u$ path has at least $L$ edges. 
Before conditioning on $E_t$, its images form a balanced walk along the prescribed colors. 
Apply~\cref{lem:initial state} averaged over the root image, showing that for every $W\subseteq V$ the probability of placing $u$ in $W$ differs from $|W\cap R_t|/m_t$ by at most $\eta$. 
Conditioning on an event whose complement has probability $\delta$ changes any event probability by at most $\delta$.  
Thus the accepted embedding satisfies
\begin{equation}\label{eq:tree-uniform}
\left|\Pr\bigl(f(u)\in W\mid\mathcal F_t\bigr)
-\frac{|W\cap R_t|}{m_t}\right|\le2\eta
\end{equation}
for every non-special $u$ and every $W\subseteq V$. 

\begin{claim}\label{claim:tree-density}
With probability $1-o(1)$, all pieces are embedded and \eqref{eq:tree-invariant} holds at every time $0\le t\le s$.
\end{claim}
\begin{proof}
Fix $W\in\mathcal N$. 
At each step, put $Y_t(W)=|W\cap R_t|/m_t$ and let $X_t(W)=|f(V(F_{t+1}))\cap W|$. 
By \eqref{eq:tree-uniform} and the bound of $B$ special vertices per piece, we have $\left|\mathbb{E}[X_t(W)\mid\mathcal F_t]-a_tY_t(W)\right|\le2\eta a_t+B$.
Also, deleting the piece gives the exact identity $Y_{t+1}(W)-Y_t(W)=(a_tY_t(W)-X_t(W))/m_{t+1}.$
Hence the total absolute conditional degree change is at most
\begin{equation}\label{eq:tree-change}
\sum_{t<s}\frac{2\eta a_t+B}{m_{t+1}}
\le\frac{2\eta n+Bs}{\ell}
\le\frac{6\eta}{\lambda}+o(1)<\frac\alpha8,
\end{equation}
where we used $\eta=\alpha\lambda/1000$.

Let $\tau$ be the first $t\le s$ for which \eqref{eq:tree-invariant} fails, with $\tau=s+1$ if it never fails. 
Write $\widetilde Y_t(W)=Y_{\min\{t,\tau\}}(W)$ and define
\[
M_t(W)=\widetilde Y_t(W)-Y_0(W)
-\sum_{i<t}\mathbb{E}[\widetilde Y_{i+1}(W)-\widetilde Y_i(W)\mid\mathcal F_i].
\]
This is a martingale. 
Moreover, $|M_{t+1}(W)-M_t(W)|\le 2a_t/m_{t+1}$ since $0\le X_t(W)\le a_t$ and $0\le Y_t(W)\le1$.
Moreover, we have $\sum_{t<s}\left(2a_t/m_{t+1}\right)^2\le4/\ell^2\sum_{t<s}a_t^2
\le 4An/\ell^2.$

If $W$ witnesses a failure of \eqref{eq:tree-invariant}, its density has fallen from at least $1/2+\alpha$ to below $1/2+\alpha/2$. 
By \eqref{eq:tree-change}, this forces $M_s(W)<-\alpha/4$. ~\cref{lem:Azuma} bounds its probability by $2\exp\left(-\frac{\alpha^2\ell^2}{128An}\right)=\exp\bigl(-\Omega(n^{2/3})\bigr)$, using $\ell\ge\lambda n/3$. 
A union bound over at most $n(n-1)$ sets $W$ proves the claim.
\end{proof}

For every deleted leaf $z\in\mathcal{L}$, let $p(z)$ be its unique neighbor in $F$, and put $c(z)=\chi(p(z)z)$. 
Define $S_v=\sum_{z\in\mathcal{L}}\mathbf 1_{\{f(p(z))v\in E(G_{c(z)})\}}$, $v\in V.$
Thus, $S_v$ counts the leaves whose prescribed parent edge allows them to be placed at $v$.

\begin{claim}\label{claim:tree-score}
With probability $1-o(1)$, the process completes and $S_v\ge(1/2+\alpha/8)\ell$ for every $v\in V$.
\end{claim}
\begin{proof}
Fix $v\in V$. 
At an active step $t$, let $Z_{t,v}=\sum_{z\in\mathcal{L},p(z)\in V(F_{t+1})}\mathbf{1}_{\{f(p(z))v\in E(G_{c(z)})\}}$ be the contribution to $S_v$ from leaves with parents in $F_{t+1}$ and set $Z_{t,v}=0$ if $t\ge\tau$. 
Let $b_t$ count those leaves whose parents in this piece are non-special. 
Applying \eqref{eq:tree-uniform} to $W=N_{G_{c(z)}}(v)$ and using \eqref{eq:tree-invariant} gives $\mathbb{E}[Z_{t,v}\mid\mathcal F_t]\ge\sum_{z\in\mathcal{L},p(z)\in V(F_{t+1}),p(z) \text{non-special}}(1/2+\alpha/2-2\eta)\ge(1/2+\alpha/2-2\eta)b_t,$ $t<\tau$.
At most $DBs\le\alpha\ell/12$ leaves have special parents, which implies that at least $\ell-DBs$ leaves have non-special parents.
Therefore, on $\{\tau>s\}$, we have
$\sum_{t<s}\mathbb{E}[Z_{t,v}\mid\mathcal F_t]\ge(1/2+\alpha/2-2\eta)(\ell-DBs)\ge(1/2+\alpha/4)\ell.$

A piece of order $a_t$ contains the parents of at most $Da_t$ deleted leaves. 
Thus $0\le Z_{t,v}\le Da_t$, including when several leaves have the same parent. 
Define $M_{j,v}=\sum_{t<j}\bigl(Z_{t,v}-\mathbb{E}[Z_{t,v}\mid\mathcal F_t]\bigr)$.
We have $M_{j+1,v}=M_{j,v}+Z_{j,v}-\mathbb E[Z_{j,v}\mid\mathcal F_j].$
Hence, $\mathbb E[M_{j+1,v}\mid\mathcal F_j]= M_{j,v}+\mathbb E[Z_{j,v}\mid\mathcal F_j]-\mathbb E\left[\mathbb E[Z_{j,v}\mid\mathcal F_j]\mid\mathcal F_j\right]=M_{j,v}$, which implies $(M_{j,v})_{j=0}^s$ is a martingale and has increments of absolute value at most $Da_t$, with $\sum_{t<s}(Da_t)^2\le D^2An$.
On $\{\tau>s\}$, we have $\sum_t Z_{t,v}=S_v$. 
Hence the event $\{\tau>s,\ S_v<(1/2+\alpha/8)\ell\}$ implies $M_{s,v}<-\alpha\ell/8$. 
~\cref{lem:Azuma} gives
\[
\Pr\bigl(\tau>s,\ S_v<(1/2+\alpha/8)\ell\bigr)
\le2\exp\left(-\frac{\alpha^2\ell^2}{128D^2An}\right)
=\exp\bigl(-\Omega(n^{1/2})\bigr).
\] 
Combine it with~\cref{claim:tree-density} and taking a union bound over $v\in V$ proves the claim.
\end{proof}

We fix an outcome satisfying both claims and let $R_s$ be the remaining set of $\ell$ host vertices. 
To embed the deleted leaves, we define a bipartite graph $\mathcal B$ with parts $\mathcal L$ and $R_s$, where $zv$ is an edge precisely when
$f(p(z))v\in E(G_{c(z)})$.
By \eqref{eq:tree-invariant}, every $z\in\mathcal{L}$ has at least $(1/2+\alpha/2)\ell$ neighbors in $R_s$, while ~\cref{claim:tree-score} gives at least $(1/2+\alpha/8)\ell$ neighbors for every $v\in R_s$.
By Hall's theorem, there exists a perfect matching in $\mathcal{B}$. 
Use this perfect matching to embed the deleted leaves. 
Leaves sharing a parent are separate vertices of $\mathcal{L}$, so they are matched to distinct host vertices, each in its own prescribed color. 
All host vertices are now used exactly once, which proves the lemma.
\end{proof}

\subsection{Paths with prescribed endpoints}\label{subsec:tree with few leaves}
We next prove ~\cref{lem:linkage}. 
The proof follows the same sequence as the Hamilton-cycle proof: choose a reservoir, set aside an absorber, find an almost cover, incorporate the remaining vertices, and finally absorb the unused reservoir. 
Here the main difficulty is to coordinate absorption with the exact allocation of internal positions, so we keep track of the internal positions on each path.

\begin{proof}[Proof of ~\cref{lem:linkage}]
We discard unused colors and vertices of $X$ that do not occur as endpoints.
We may assume $0<\alpha<1/4$.

Take $q=q(2,\alpha)$ and choose $\theta$ as in~\cref{lem:block absorb}. 
Choose constants
\[
0<1/N\ll1/K\ll1/L\ll\varepsilon\ll\gamma\ll\beta\ll\theta,1/q,\alpha.
\]
The short paths in the almost-cover step will have $L$ vertices.
Put $b=\lfloor\beta N\rfloor,\ g=\lfloor\gamma N\rfloor,
\ T=\lceil b/\theta\rceil+b$ and $z=b+g$. 
Since $\sum_{j=1}^{p}\ell_j=N$ and $\ell_j\ge K$, we have $p\le N/K$. After discarding unused colors and endpoints, we may therefore assume
that $m\le N+p\le 2N,\ |X|\le 2p\le 2N/K$.
Consider the family $\mathcal N=\bigl\{N_{G_i}(v)\cap U:i\in[m],\ v\in U\cup X\bigr\}.$
By assumption, every $Q\in\mathcal N$ satisfies $|Q|\ge (1/2+\alpha)N$,
and $|\mathcal N|\le m(N+|X|)\le 4N^2$.
Choose a $z$-subset $Z$ of $U$ uniformly at random.
For every fixed $Q\in\mathcal N$, the random variable
$X:=|Q\cap Z|$ is hypergeometrically distributed
with mean $\mu=z|Q|/N\ge (1/2+\alpha)z.$
Applying~\cref{lem:chernoff} with $\varepsilon=\alpha z/(2\mu)$ gives
$\Pr\left(|Q\cap Z|<(1/2+\alpha/2)z\right)\le 2\exp\left(-\frac{\alpha^2z^2}{12\mu}\right)\le 2\exp\left(-\frac{\alpha^2z}{12}\right).$
Since $z=\Theta(N)$, a union bound over $\mathcal N$ shows that, for sufficiently large $N$, there exists a $z$-set $Z\subseteq U$ satisfying $\delta(\mathcal G,Z)\ge (1/2+\alpha/2)|Z|$.
Consequently, for any $u,v\in U\cup X$ and any $i,j\in[m]$, we have
$|N_{G_i}(u)\cap N_{G_j}(v)\cap Z|\ge d_{G_i}(u,Z)+d_{G_j}(v,Z)-|Z|\ge \alpha|Z|.$
We use this observation repeatedly to make two-edge connections.

We first prepare the absorbing paths.
Recall that each $q$-block uses $q+1$ colors, and we reserve two further colors to connect it to the preceding part of its path.
Thus path $j$ can accommodate $c_j=\left\lfloor(\ell_j-1)/(q+3)\right\rfloor$ blocks while leaving at least two colors to reach $y_j$.
Since $\ell_j\ge K\gg q$, we have $(q+3)c_j>\ell_j-q-4\ge\ell_j/2,$ and hence $\sum_{j\in[p]}c_j>T$.
We assign the $T$ blocks in order to the paths, filling each path to its capacity before moving to the next.
Let $J\subseteq[p]$ denote the set of paths receiving blocks, and let $t_j\ge1$ be the number assigned to path $j\in J$. 
All but the last such path are filled to capacity.
We denote the index of the last path as $j'\in J$, hence for every $j\in J\setminus\{j'\}$, we have $c_j=t_j$ and $\sum_{j\in J\setminus\{j'\}}t_j\le T$.
Since the last selected path has length at most $2K$, it follows that $\sigma:=\sum_{j\in J}\ell_j\le2(q+3)\sum_{j\in J\setminus\{j'\}}t_j+2K\le 2(q+3)T+2K<\alpha N/10$.

For each $j\in J$, we reserve a two-edge connector starting at $x_j$ before each of its $t_j$ blocks. 
The final $\ell_j+1-t_j(q+3)\ge2$ edges are reserved for the connection to $y_j$.
This determines consecutive, disjoint color intervals for all the blocks and connectors.
We choose pairwise distinct endpoints for the $T$ blocks in $U\setminus Z$ and let $D$ be their union. 
Apply ~\cref{lem:block absorb} to $(G_i[U])_i$, these blocks and $Z$. 
It gives $A_0\subseteq U\setminus(Z\cup D)$ of size $qT-b$ such that for every $Z'\in\binom Zb$, all blocks admit disjoint internal solutions whose union is $A_0\cup Z'$.

We now fix the connectors on the paths in $J$, taking all their internal vertices in $U\setminus(A_0\cup Z\cup D)$ and avoiding those already chosen. A prescribed connector of length at least two can be built by choosing its initial internal vertices greedily, then choosing its last internal vertex in the common neighborhood of the current vertex and the prescribed endpoint in the last two colors. 
This common neighborhood in $U$ has size at least $2\alpha N$, whereas the total number of excluded vertices is at most $|A_0|+|Z|+|D|+\sigma=O_{q,\theta}(\beta N)+O(K)\ll\alpha N$. 
Thus all connectors can be fixed disjointly.
Let $A$ consist of $A_0$, the block endpoints and these connector vertices. For any $Z'\in\binom Zb$, the blocks use precisely $A_0\cup Z'$, while all other internal vertices are fixed. 
Thus $|A|=\sigma-b$, and the prescribed paths in $J$ have pairwise disjoint internal vertex sets partitioning $A\cup Z'$.
The arrangement of the $q$-blocks, two-edge connectors, and final tail along each prescribed path in $J$ is illustrated in Figure~\ref{fig:block-path}.

\begin{figure}[htbp]
    \centering
    \includegraphics[width=\linewidth]{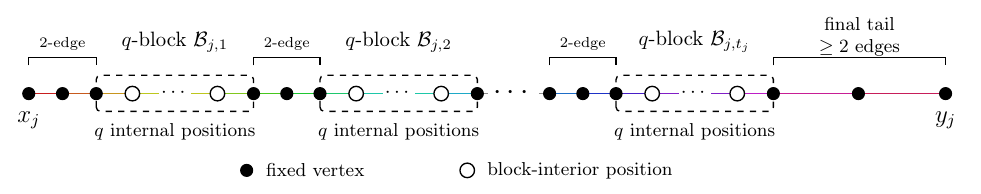}
    \caption{Schematic of the absorbing construction
    on a selected path $j\in J$.}
    \label{fig:block-path}
\end{figure}

It remains to construct the paths indexed by $I:=[p]\setminus J$. 
Write $M:=\sum_{j\in I}\ell_j=N-\sigma$ and $W:=U\setminus(A\cup Z)$, so $n':=|W|=M-g\ge N/2$. 
As $|U\setminus W|=\sigma+g<\alpha N/5$, every relevant graph induced on $W$ has minimum degree at least $(1/2+\alpha/2)n'$. 
Let $s=\lfloor(1-\varepsilon)n'/L\rfloor$ and $d=n'-sL\le\varepsilon N+L$. Choose integers $s_j$ for $j\in I$ with $\sum_{j\in I}s_j=s$ and $|s_j-s\ell_j/M|\le1$. 
Our choice of $K$ ensures $s_j\ge1$ for every $j\in I$.

On path $j\in I$, reserve $s_j$ short paths of order $L$, each preceded by one vertex of $Z$ which gives a two-edge connection from $x_j$ or the previous short path. 
This uses $(L+1)s_j$ internal positions, leaving a tail of $r_j:=\ell_j-(L+1)s_j$ internal vertices before $y_j$. 
Since $s\le(M-g)/L$, we have $r_j\ge\ell_j-(L+1)(1+s\ell_j/M)\ge\ell_j(g/M-1/L)-(L+1)\ge1$, using $g/M\ge\gamma/2$, $1/L\ll\gamma$ and $K\gg L/\gamma$. 
The positions of the short paths now specify their color lists. 
At most $s(L-1)<n'$ distinct colors occur in these lists. 
Relabeling them as $[m']$ with $m'\le n'$, we may apply~\cref{lem:almost cover} with $k=2$ to obtain $s$ pairwise disjoint, exactly colored $L$-vertex paths in $W$. 
We join these short paths in their prescribed order, starting from the corresponding $x_j$ and using one unused vertex of $Z$ before each short path.
This uses $s$ vertices of $Z$.

There remain exactly $d$ unused vertices of $W$. 
We assign them to the tails, allowing two internal positions for each such vertex and one final position to reach $y_j$. 
Indeed, if $p'=|I|$, then $\sum_{j\in I}r_j=M-(L+1)s=d+g-s$, and hence
\[
\sum_{j\in I}\left\lfloor\frac{r_j-1}{2}\right\rfloor
\ge\frac{d+g-s-2p'}{2}\ge d.
\]
The last inequality follows from $d+s+2p'\le\varepsilon N+L+N/L+2N/K\le g$. We can therefore distribute the $d$ vertices among the tails, giving $d_j$ vertices to tail $j$ so that $\sum_{j\in I}d_j=d$ and $2d_j+1\le r_j$.

We complete tail $j$ by starting at its last short path and visiting the $d_j$ assigned vertices of $W$ in any order.
Each visit uses one fresh vertex of $Z$ for each two-edge connection. 
We then extend through $r_j-2d_j-1$ further vertices of $Z$ and use one last vertex of $Z$ to make a two-edge connection to $y_j$. 
Throughout this construction, all edges are taken in the prescribed colors. 
The completed tail has exactly $r_j$ internal vertices, of which $d_j$ belong to $W$. 
These choices are always possible: altogether the paths in $I$ use only
$s+\sum_{j\in I}(r_j-d_j)=s+(d+g-s)-d=g\ll\alpha|Z|$
vertices of $Z$, and each required one- or two-edge extension has more available choices than previously used reservoir vertices.

We have now covered $W$ and exactly $g$ vertices of $Z$ with the paths in $I$. 
The remaining set $Z'\subseteq Z$ has size $b$, so the absorbing property completes all paths in $J$ on $A\cup Z'$. 
Their internal vertex sets, together with those of the paths in $I$, partition $U$. 
Every path has its prescribed endpoints, length and color order, as required.
\end{proof}

\section{A spanning-tree construction}\label{sec:construction}

In this section, we prove~\cref{thm:tree-construction}.

\begin{proof}[Proof of~\cref{thm:tree-construction}]
    Fix $\ell\ge 9$. 
    Choose a positive integer $q$ divisible by $3\ell(24\ell+1)$, and set $n:=q^2,\ h:=\frac{n}{24\ell+1},\ B:=4\ell h.$ 
    Then $h$ and $B$ are integers, $3\ell\mid q$, and $n=6B+h.$

    We first describe the common vertex set. 
    Partition $ V=A\mathbin{\dot\cup}P_1\mathbin{\dot\cup}P_2 \mathbin{\dot\cup}P_3\mathbin{\dot\cup}P_4$ such that $|A|=2B+h$ and $|P_i|=B$ for every $i\in[4]$.
    For each $i\in[4]$, further partition $P_i=X_i\mathbin{\dot\cup}R_i(1)\mathbin{\dot\cup}R_i(2)\mathbin{\dot\cup}R_i(3),$ where $|X_i|=|R_i(c)|=\ell h=B/4$ for every $c\in[3]$.
    Finally, partition $R_i(c)=\mathbin{\dot\bigcup}_{t=1}^{\ell}R_i(c,t),\ |R_i(c,t)|=h$, and write $R_i:=R_i(1)\cup R_i(2)\cup R_i(3).$
    Thus $|R_i|=\frac{3B}{4}.$
    Let $\mathcal C:=\{0\}\mathbin{\dot\cup}([3]\times[\ell]).$
    
    We now define the graph collection.
    Set $ G_0:=K_{A,V\setminus A}\cup\bigcup_{i=1}^4K_{P_i}.$
    For $(c,t)\in[3]\times[\ell]$, define $ G_{c,t}:=K_A\cup\bigcup_{i=1}^4K_{R_i}\cup\bigcup_{1\le i<j\le4}K_{P_i,P_j}\cup\bigcup_{i=1}^4K_{A,R_i(c)}\cup\bigcup_{i=1}^4K_{X_i,R_i(c,t)}.$
    No other edges are present.

    We include the two properties of the construction which will be used later. 
    If $v\in R_i(c',t'),$ then  
    \[
    N_{G_{c,t}}(v)\cap A=
    \begin{cases}
        A, & c=c',\\
        \varnothing, & c\ne c',
    \end{cases}
    \]
    and 
    \[
    N_{G_{c,t}}(v)\cap X_i=
    \begin{cases}
        X_i, & (c,t)=(c',t'),\\
        \varnothing, & (c,t)\ne(c',t').
    \end{cases}
    \]
    Moreover, no graph $G_{c,t}$ contains an edge between $A$ and $X_i$, or an edge inside $X_i$.

    We next verify the minimum degree. 
    In $G_0$, every vertex of $A$ has degree $4B$, while every vertex of $P_i$ has degree $|A|+|P_i|-1=3B+h-1.$
    Hence $\delta(G_0)=3B+h-1.$
    Fix $(c,t)\in[3]\times[\ell]$. A vertex of $A$ has degree $|A|-1+\sum_{i=1}^4|R_i(c)|=3B+h-1.$
    A vertex of $X_i$ has degree $3B+h$, while every vertex of $R_i$ has degree at least $3B+3B/4-1.$ 
    Since $ 3B/4=3\ell h\ge h,$ we obtain $\delta(\mathcal G)=3B+h-1=(n+h)/2-1=\left(\frac12+\frac{1}{2(24\ell+1)}\right)n-1.$

    We now define the tree and its prescribed color pattern. 
    The tree $T$ consists of a root $r$ with $q-1$ neighbours, called centres, each of which has $q$ additional leaf neighbours. 
    Then $|V(T)|=1+(q-1)+q(q-1)=q^2=n$ and $\Delta(T)=q+1=\sqrt n+1.$

    Every root--centre edge has color $0$.
    At each centre, exactly $q/(3\ell)$ leaf edges receive each color $(c,t)\in[3]\times[\ell].$
    This is possible since $3\ell\mid q$, and defines the prescribed pattern $\chi:E(T)\to \mathcal{C}$.

    Suppose, for a contradiction, that there exists an exactly $\chi$-colored copy of $T$ in $\mathcal G$. 
    Let $ f:V(T)\to V$ be the corresponding bijection. 
    We distinguish two cases according to the image of the root.
    
    \noindent
    \textbf{Case 1.} $f(r)\in A$.
    
    Since every root--centre edge has color $0$, the definition of $G_0$ implies that every centre is mapped into $V\setminus A$.
    
    A centre mapped into $X_i$ has no leaf-edge neighbour in $A$. If a centre is mapped into $R_i(c',t')$, then only the
    $\ell$ colors $(c',1),\ldots,(c',\ell)$ may be used on leaf edges entering $A$. 
    Hence each centre sends at most $\ell\cdot q/(3\ell)=q/3$ leaves into $A$. 
    Therefore $|f(V(T))\cap A|\le1+(q-1) q/3<q^2/3=n/3<|A|.$
    This contradicts the bijectivity of $f$.
    
    \noindent
    \textbf{Case 2.} $f(r)\in P_i$ for some $i\in[4]$.
    
    By the definition of $G_0$, every centre is mapped into $A\cup P_i$.
    We estimate the number of vertices which can be mapped into $X_i$.
    
    A centre mapped into $A$ or $X_i$ cannot send a leaf into $X_i$.
    If a centre is mapped into $R_i(c',t')$, then only the
    single color $(c',t')$ may be used on a leaf edge entering $X_i$.
    Thus each centre sends at most $q/(3\ell)$ leaves into $X_i$.
    
    The root and the centres together account for only $q$ vertices.
    Consequently, $|f(V(T))\cap X_i|\le q+(q-1)\frac{q}{3\ell}.$
    
    Since $\ell\ge9$ and $q\ge3\ell(24\ell+1)$,
    \[
    \left(\frac{\ell}{24\ell+1}-\frac{1}{3\ell}\right)q=\frac{(3\ell^2-24\ell-1)q}{3\ell(24\ell+1)}>1.
    \]
    Hence $q+(q-1)\frac{q}{3\ell}<q+\frac{q^2}{3\ell} <\frac{\ell q^2}{24\ell+1}=\ell h=|X_i|.$
    Again this contradicts the bijectivity of $f$.
    
    The two cases exhaust all possible images of the root. 
    Therefore $\mathcal G$ contains no exactly $\chi$-colored copy of $T$.
    Since $q$ can be chosen to be an arbitrarily large multiple of
    $3\ell(24\ell+1)$, the proof is complete.
\end{proof}

\section{Concluding remarks}\label{sec:conluding rmk}
Our results add two further examples to the study of transversal and pattern thresholds. 
For tight Hamilton cycles in uniform hypergraphs, prescribing the entire color pattern does not change the asymptotic minimum codegree threshold. However, for spanning trees, our construction gives a positive linear separation even for trees of maximum degree $\sqrt n+1$ and a fixed palette of $28$ colors.

The proofs and the construction suggest two different mechanisms.
For tight Hamilton cycles, we can connect ordered ends through prescribed sequences of colors and exchange vertices inside absorbing blocks without changing either their endpoints or their color lists. 
In the tree construction, the colors incident with a branching vertex restrict the parts of the host in which its neighbors can be embedded. 
Summing these restrictions over many branching vertices yields a global obstruction to a spanning embedding.
Thus the structure of the prescribed pattern, as well as the maximum degree of the target, may be relevant to pattern universality.

It would be interesting to determine the exact pattern threshold for tight Hamilton cycles.
In particular, for every fixed $k\ge3$, does there exist a constant $C_k$ such that $\delta_{k-1}(\mathcal H)\ge n/2+C_k$ guarantees every prescribed pattern of a tight Hamilton cycle for all sufficiently large $n$?
For spanning trees, the main gap lies in the permitted maximum degree.
\Cref{thm:tree} allows $\Delta(T)\le n^{c(\alpha)}$, whereas
\cref{thm:tree-construction} gives an obstruction at
$\Delta(T)=\sqrt n+1$ for a fixed positive degree margin.
It would be particularly interesting to determine whether pattern
universality still holds under $\delta(\mathcal G)\ge(1/2+\alpha)n$ for all trees with $\Delta(T)=o(\sqrt n)$, or whether the separation can already occur at a smaller polynomial maximum degree. 

Finally, our approach may be useful beyond the structures considered here.
The balanced-extension walk preserves the uniform distribution on ordered
$(k-1)$-tuples even when successive edges are drawn from different members
of a collection; its coupling property gives connections respecting an
arbitrary prescribed color sequence. 
Together with fixed-pattern swap absorption, this suggests a possible approach to other spanning problems, such as Hamilton $\ell$-cycles and prescribed-endpoint systems of tight paths.
The different overlap and divisibility constraints in these settings would
require corresponding connecting and absorbing structures. 
Understanding these structures may help clarify when transversal and pattern thresholds coincide.

\medskip
\noindent
\subsection*{Acknowledgement and statement on AI use.}
The author wishes to thank Jie Han for valuable advice on the presentation of the paper.

GPT 6-Pro was used as a supporting tool to suggest a candidate construction and to assist in the development of parts of a proof based on the ideas and guidance of the author. 
Furthermore, GPT is also utilized to enhance the English expression of this article.
All mathematical arguments were independently verified and refined by the author, who takes full responsibility for the manuscript.

\bibliographystyle{abbrv}
\bibliography{ref}

@article{glock2021counting,
  title={Counting Hamilton cycles in Dirac hypergraphs},
  author={Glock, Stefan and Gould, Stephen and Joos, Felix and K{\"u}hn, Daniela and Osthus, Deryk},
  journal={Combinatorics, Probability and Computing},
  volume={30},
  number={4},
  pages={631--653},
  year={2021},
  publisher={Cambridge University Press}
}

@book{frieze2015introduction,
  title={Introduction to random graphs},
  author={Frieze, Alan and Karo{\'n}ski, Micha{\l}},
  volume={1},
  year={2015},
  publisher={Cambridge University Press Cambridge}
}

@book{AlonSpencer2016,
  author    = {Noga Alon and Joel H. Spencer},
  title     = {The Probabilistic Method},
  edition   = {4},
  publisher = {John Wiley \& Sons},
  address   = {Hoboken, NJ},
  year      = {2016}
}

@article{reiher2019minimum,
  title={Minimum vertex degree condition for tight Hamiltonian cycles in 3-uniform hypergraphs},
  author={Reiher, Christian and R{\"o}dl, Vojt{\v{e}}ch and Ruci{\'n}ski, Andrzej and Schacht, Mathias and Szemer{\'e}di, Endre},
  journal={Proceedings of the London Mathematical Society},
  volume={119},
  number={2},
  pages={409--439},
  year={2019},
  publisher={Wiley Online Library}
}

@article{lo2015f,
  title={F-factors in hypergraphs via absorption},
  author={Lo, Allan and Markstr{\"o}m, Klas},
  journal={Graphs and Combinatorics},
  volume={31},
  number={3},
  pages={679--712},
  year={2015},
  publisher={Springer}
}

@article {MR47308,
    AUTHOR = {Dirac, G. A.},
     TITLE = {Some theorems on abstract graphs},
   JOURNAL = {Proc. London Math. Soc. (3)},
  FJOURNAL = {Proceedings of the London Mathematical Society. Third Series},
    VOLUME = {2},
      YEAR = {1952},
     PAGES = {69--81},
      ISSN = {0024-6115,1460-244X},
   MRCLASS = {56.0X},
  MRNUMBER = {47308},
MRREVIEWER = {W.\ T.\ Tutte},
       DOI = {10.1112/plms/s3-2.1.69},
       URL = {https://doi.org/10.1112/plms/s3-2.1.69},
}

@article{rodl2008approximate,
  title={An approximate Dirac-type theorem for k-uniform hypergraphs},
  author={R{\"o}dl, Vojt{\v{e}}ch and Szemer{\'e}di, Endre and Ruci{\'n}ski, Andrzej},
  journal={Combinatorica},
  volume={28},
  number={2},
  pages={229--260},
  year={2008},
  publisher={Springer}
}

@article{joos2020rainbow,
  title={On a rainbow version of Dirac's theorem},
  author={Joos, Felix and Kim, Jaehoon},
  journal={Bulletin of the London Mathematical Society},
  volume={52},
  number={3},
  pages={498--504},
  year={2020},
  publisher={Wiley Online Library}
}

@article{cheng2025transversal,
  title={Transversal Hamilton cycle in hypergraph systems},
  author={Cheng, Yangyang and Han, Jie and Wang, Bin and Wang, Guanghui and Yang, Donglei},
  journal={SIAM Journal on Discrete Mathematics},
  volume={39},
  number={1},
  pages={55--74},
  year={2025},
  publisher={SIAM}
}

@inproceedings{cheng2023rainbow,
  title={Rainbow spanning structures in graph and hypergraph systems},
  author={Cheng, Yangyang and Han, Jie and Wang, Bin and Wang, Guanghui},
  booktitle={Forum of Mathematics, Sigma},
  volume={11},
  pages={e95},
  year={2023},
  organization={Cambridge University Press}
}

@article{gupta2023general,
  title={A general approach to transversal versions of Dirac-type theorems},
  author={Gupta, Pranshu and Hamann, Fabian and M{\"u}yesser, Alp and Parczyk, Olaf and Sgueglia, Amedeo},
  journal={Bulletin of the London Mathematical Society},
  volume={55},
  number={6},
  pages={2817--2839},
  year={2023},
  publisher={Wiley Online Library}
}

@article{bowtell2025universality,
  title={Universality for transversal Hamilton cycles},
  author={Bowtell, Candida and Morris, Patrick and Pehova, Yanitsa and Staden, Katherine},
  journal={Bulletin of the London Mathematical Society},
  volume={57},
  number={3},
  pages={711--729},
  year={2025},
  publisher={Wiley Online Library}
}

@article{chakraborti2023bandwidth,
  title={A bandwidth theorem for graph transversals},
  author={Chakraborti, Debsoumya and Im, Seonghyuk and Kim, Jaehoon and Liu, Hong},
  journal={arXiv preprint arXiv:2302.09637},
  year={2023}
}

@article{chakraborti2024hamilton,
  title={Hamilton transversals in tournaments},
  author={Chakraborti, Debsoumya and Kim, Jaehoon and Lee, Hyunwoo and Seo, Jaehyeon},
  journal={Combinatorica},
  volume={44},
  number={6},
  pages={1381--1400},
  year={2024},
  publisher={Springer}
}

@article {MR4451150,
    AUTHOR = {Montgomery, Richard and M\"uyesser, Alp and Pehova, Yani},
     TITLE = {Transversal factors and spanning trees},
   JOURNAL = {Adv. Comb.},
  FJOURNAL = {Advances in Combinatorics},
      YEAR = {2022},
     PAGES = {Paper No. 3, 25},
      ISSN = {2517-5599},
   MRCLASS = {05C35},
  MRNUMBER = {4451150},
       DOI = {10.19086/aic.2022.3},
       URL = {https://doi.org/10.19086/aic.2022.3},
}

@article{heath2025universality,
  title={Universality for transversal powers of Hamilton cycles},
  author={Heath, Emily and Hyde, Joseph and Morrison, Natasha and Ogden, Shannon},
  journal={arXiv preprint arXiv:2510.18163},
  year={2025}
}

@article{christoph2025universality,
  title={Universality for transversal Hamilton cycles in random graphs},
  author={Christoph, Micha and Martinsson, Anders and Milojevi{\'c}, Aleksa},
  journal={arXiv preprint arXiv:2505.05385},
  year={2025}
}

@article{cuckler2009hamiltonian,
  title={Hamiltonian cycles in Dirac graphs},
  author={Cuckler, Bill and Kahn, Jeff},
  journal={Combinatorica},
  volume={29},
  number={3},
  pages={299--326},
  year={2009},
  publisher={Springer}
}

@book{levin2026markov,
  title={Markov chains and mixing times},
  author={Levin, David A and Peres, Yuval},
  year={2017},
  publisher={American Mathematical Society}
}

@article{sudakov2010randomized,
  title={A randomized embedding algorithm for trees},
  author={Sudakov, Benny and Vondr{\'a}k, Jan},
  journal={Combinatorica},
  volume={30},
  number={4},
  pages={445--470},
  year={2010},
  publisher={Springer}
}

@article{montgomery2020embedding,
  title={Embedding rainbow trees with applications to graph labelling and decomposition},
  author={Montgomery, Richard and Pokrovskiy, Alexey and Sudakov, Benjamin},
  journal={J. Eur. Math. Soc},
  volume={22},
  number={10},
  pages={3101--3132},
  year={2020}
}

\end{document}